\documentclass[12pt,twoside]{amsart}
\usepackage{amssymb,amsthm,amsmath}
\usepackage[ngerman, english]{babel}
\usepackage{amscd}
\usepackage{hyperref}
\usepackage{color}

\usepackage[normalem]{ulem} 
\newcommand\redout{\bgroup\markoverwith
{\textcolor{red}{\rule[.5ex]{2pt}{0.4pt}}}\ULon}

\title[Zariskian adic spaces]
{Zariskian adic spaces} 
\author{Hiromu Tanaka} 
\address{Graduate School of Mathematical Sciences, 
the University of Tokyo, 3-8-1 Komaba, Meguro-ku, Tokyo 153-8914, Japan.} 
\email{tanaka@ms.u-tokyo.ac.jp}

\newcommand{\hd}[0]{{\operatorname{hd}}}

\newcommand{\Zar}[0]{{\operatorname{Zar}}}

\newcommand{\Ker}[0]{{\operatorname{Ker}}}

\newcommand{\Image}[0]{{\operatorname{Im}}}

\newcommand{\Spa}[0]{{\operatorname{Spa}}}

\newcommand{\Proj}[0]{{\operatorname{Proj}}}

\newcommand{\Spec}[0]{{\operatorname{Spec}}}

\newtheorem{thm}{Theorem}[section]
\newtheorem{lem}[thm]{Lemma}
\newtheorem{cor}[thm]{Corollary}
\newtheorem{prop}[thm]{Proposition}

\newtheorem{step}{Step}
\newtheorem*{claim}{Claim}

\theoremstyle{definition}

\newtheorem{ex}[thm]{Example}

\newtheorem{dfn}[thm]{Definition}

\newtheorem{rem}[thm]{Remark}
      
\newtheorem{nota}[thm]{Notation}

\makeatletter
  
  \@addtoreset{equation}{section}
  \makeatother

\newcommand{\MO}{\mathcal{O}}
\newcommand{\F}{\mathbb{F}}

\newcommand{\R}{\mathbb{R}}
\newcommand{\Q}{\mathbb{Q}}
\newcommand{\Z}{\mathbb{Z}}
\newcommand{\p}{\mathfrak{p}}

\newcommand{\m}{\mathfrak{m}}

\begin{document}

\maketitle

\begin{abstract}
We introduce a Zariskian analogue of the theory of Huber's adic spaces. 
\end{abstract}


\tableofcontents

\section{Introduction}

Tate introduced a $p$-adic analytic geometry 
so-called the rigid geometry. 
In the original definition by Tate, 
a rigid analytic space is not a topological space but 
a Grothendieck topological space (cf. \cite{BGR84}). 
To remedy this situation, 
Huber established the theory of adic spaces. 
He introduced a topological space $\Spa\,A$, 
called an affinoid spectrum, associated to an affinoid ring $A=(A^{\rhd}, A^+)$, 
where $A^{\rhd}$ is an f-adic ring and 
$A^+$ is a certain open subring of $A^{\rhd}$. 
An adic space is obtained by gluing affinoid spectra. 

Although Huber introduced a structure presheaf $\MO_A$ 
on $\Spa\,A$ for an arbitrary affinoid ring $A$, 
it is not a sheaf in general 
(cf. \cite{BV}, \cite[the example after Proposition 1.6]{Hub94}, \cite{Mih16}). 
The reason for this is that 
the definition of $\MO_A$ depends on the completion, 
which is a transcendental operation. 
Thus it is natural to ask whether 
the theory can be more well-behaved 
after replacing the completion by a more algebraic operation. 
For example, if $(A, I)$ is a pair consisting of a ring $A$ and an ideal $I$ of $A$, 
then we can associate the Zariskian ring $A^{\Zar}$ and 
the henselisation $A^h$ with respect to $I$. 
The henselisation $A^h$ is known as an algebraic approximation of the $I$-adic completion $\widehat{A}$, 
whilst the associated Zariskian ring $A^{\Zar}$ is closer to the original ring $A$ than $A^h$. 
The purpose of this paper is to establish the Zariskian version 
of Huber's theory. 
Therefore, the first step is to introduce a notion of Zariskian f-adic rings. 

\begin{dfn}[Definition~\ref{d-zar}, Remark \ref{r-zar1}, Definition~\ref{d-zar-top}]\label{intro-d-zar}
Let $A$ be an f-adic ring. 
\begin{enumerate}
\item 
We set $S_A^{\Zar}:=1+A^{\circ\circ}$. 
It is easy to show that $S_A^{\Zar}$ is a multiplicative subset of $A$. 
We set $A^{\Zar}:=(S_A^{\Zar})^{-1}A$. 
Both $A^{\Zar}$ and the natural ring homomorphism $\alpha:A \to A^{\Zar}$ 
are called the {\em Zariskisation} of $A$. 
We say that $A$ is {\em Zariskian} if $\alpha:A \to A^{\Zar}$ 
is bijective. 
\item 
For a ring of definition $A_0$ of $A$ and 
an ideal of definition $I_0$ of $A_0$, 
we equip $A^{\Zar}$ with the group topology 
defined by the images of $\{I_0^kA_0^{\Zar}\}_{k\in \Z_{>0}}$. 
We can show that 
this topology does not depend on the choice of $A_0$ and $I_0$ 
(cf. Lemma~\ref{l-top-indep}). 
\end{enumerate}
\end{dfn}

We will prove that $A^{\Zar}$ satisfies some reasonable properties. 
For instance, $A^{\Zar}$ is a Zariskian f-adic ring (Theorem~\ref{t-zar-zar}) and 
$A^{\Zar}$ has the same completion as the one of $A$ 
(Theorem~\ref{t-comp-factor}). 
However one might consider that 
the definition of the topology of $A^{\Zar}$ 
is somewhat artificial. 
The following theorem asserts that 
our definition given above can be characterised 
in a category-theoretic way, 
i.e. $A^{\Zar}$ is an initial object of 
the category of Zariskian f-adic $A$-algebras.

\begin{thm}[Theorem~\ref{t-zar-univ}]\label{intro-t-zar-univ}
Let $A$ be an f-adic ring and 
let $\alpha:A \to A^{\Zar}$ be the Zariskisation of $A$. 
Then, for any continuous ring homomorphism $\varphi:A \to B$ 
to a Zariskian f-adic ring $B$, 
there exists a unique continuous ring homomorphism 
$\psi:A^{\Zar} \to B$ such that $\varphi=\alpha \circ \psi$. 
\end{thm}

For an affinoid ring $A=(A^{\rhd}, A^+)$, 
we introduce a presheaf $\MO_A^{\Zar}$ on $\Spa\,A$ 
in the same way as in Huber's theory. 
The presheaf $\MO_A$ introduced by Huber is not a sheaf in general, 
whilst the Zariskian version $\MO^{\Zar}_A$ 
is always a sheaf. 

\begin{thm}[Theorem~\ref{t-sheafy}]
For an affinoid ring $A=(A^{\rhd}, A^+)$, the presheaf $\MO_A^{\Zar}$ on $\Spa\,A$ is a sheaf. 
\end{thm}

Then one might be tempted to hope 
the Tate acyclicity in general. 
Unfortunately this is not the case.

\begin{thm}[Theorem~\ref{t-non-TA}]
There exists an affinoid ring $A=(A^{\rhd}, A^+)$ such that 
$H^1(\Spa\,A, \MO_A^{\Zar}) \neq 0$. 
\end{thm}

Although the Zariskian structure sheaf $\MO_A^{\Zar}$ 
does not behave nicely 
to establish a theory of coherent sheaves, 
the Zariskian rings might be still useful. 
For instance, if $A$ is a noetherian ring equipped with an $\m$-adic topology for some maximal ideal $\m$, 
then the Zariskisation $A^{\Zar}$ is nothing but 
the local ring $A_{\m}$ at $\m$. 
Therefore, in this situation, 
$A$ is Zariskian if and only if its $\m$-adic completion 
$A \to \widehat{A}$ is faithfully flat. 
The flatness of completion is a thorny problem 
for non-noetherian rings, 
whilst we prove that the completion is actually faithfully flat, 
under the assumption that $A$ is Zariskian and the completion is flat. 
More generally, we obtain the following result.

\begin{thm}[Corollary~\ref{c-ff-criterion}]\label{intro-t-ff}
Let $\varphi:A \to B$ be a continuous ring homomorphism 
of Zariskian f-adic rings. 
Assume that the induced map $\widehat{\varphi}:\widehat{A} \to \widehat{B}$ 
is an isomorphism of topological rings. 
Then the following hold. 
\begin{enumerate}
\item 
Any maximal ideal of $A$ is contained in the image of 
the induced map $\Spec\,B \to \Spec\,A$. 
\item 
If $\varphi$ is flat, then $\varphi$ is faithfully flat. 
\end{enumerate}
\end{thm}

Theorem~\ref{intro-t-ff} is a consequence of 
the following characterisation of Zariskian f-adic rings.

\begin{thm}[Theorem~\ref{t-characterise}]
Let $(A, A^+)$ be an affinoid ring. 
Then $A$ is Zariskian if and only if 
any maximal ideal of $A$ is contained in the image of 
the natural map 
$$\theta:\Spa\,(A, A^+) \to \Spec\,A, \quad v \mapsto \Ker(v).$$
\end{thm}

For a Zariskian affinoid ring $(A, A^+)$, 
the above theorem claims that 
the image of $\theta$ contains all the maximal ideals, 
however the map $\theta$ 
is not surjective in general (Theorem~\ref{t-non-surje}).

\medskip

\textbf{Acknowledgement:} 
The author was funded by EPSRC. 
He would like to thank the referee for reading the paper carefully and for giving many constructive comments. 

\section{Preliminaries}

In this section, we summarise notation and basic results.

\subsection{Notation}

\begin{enumerate}
\item 
Throughout the paper, a {\em ring} is 
always assumed to be commutative and to have a unity element 
with respect to the multiplication. 
We say that $P$ is a {\em pseudo-ring}, 
if $P$ satisfies all the axioms of rings except for the existence of a multiplicative identity. 
A subset $Q$ of a pseudo-ring is a {\em pseudo-subring} of $P$
if $Q$ is an additive subgroup of $P$ such that 
$q_1q_2 \in Q$ for all $q_1, q_2 \in Q$. 
In this case, we consider $Q$ as a pseudo-ring. 
For example, any ideal of a ring $A$ is a pseudo-subring of $A$. 
\item 
We will freely use the notation and terminology of 
\cite{Hub93} and \cite{Hub94}. 
In particular, given a topological ring $A$, 
we denote by 
$A^{\circ}$ (resp. $A^{\circ\circ}$) the subset of $A$ consisting of 
the power-bounded (resp. topological nilpotent) elements. 
If $A$ is an f-adic ring, then $A^{\circ}$ is an open subring of $A$ 
and $A^{\circ\circ}$ is an open ideal of $A^{\circ}$. 
A map $\varphi:A \to B$ of topological rings 
is an {\em isomorphism of topological rings} 
if $\varphi$ is bijective and both $\varphi$ and $\varphi^{-1}$ 
are continuous ring homomorphisms. 
\item 
For a partially ordered set $I$ and an element $i_0 \in I$, 
we set $I_{\geq i_0}:=\{i \in I\,|\, i\geq i_0\}$. 
We define $I_{>i_0}$ in the same way. 
\item 
For a topological space $X$, 
the topology on $X$ is {\em trivial} 
if any open subset of $X$ is equal to $X$ or the empty set. 
\end{enumerate}

\begin{dfn}\label{d-OBI}
Let $A$ and $B$ be topological rings. 
\begin{enumerate}
\item 
We define $\mathfrak O_A$ as 
the set of the open subrings of $A$. 
If $f:A \to B$ is a continuous ring homomorphism, 
we define $f^*:\mathfrak O_{B} \to \mathfrak O_{A}$ 
by $B_0 \mapsto f^*(B_0):=f^{-1}(B_0)$. 
If $f^*$ is bijective, then the inverse map is denoted by $f_*$. 
\item 
Let $\mathfrak B_{A}$ be the set of the bounded open subrings 
of $A$. 
\item 
Let $\mathfrak I_{A}$ be the set of 
the open subrings of $A$ 
that are contained in $A^{\circ}$ and integrally closed in $A$. 
\end{enumerate}
\end{dfn}

\subsection{Group topologies}

Let $P$ be a pseudo-ring. 
Let $\{P_n\}_{n \in \Z_{\geq 0}}$ be a set of additive subgroups of $P$ 
such that $P_{0} \supset P_{1} \supset \cdots$. 
Then we call the {\em group topology} on $P$ 
induced by $\{P_n\}_{n \in \Z_{\geq 0}}$ is a topology on $P$ 
such that given a subset $U$ of $P$, 
$U$ is an open subset of $P$ if and only if 
for any $x \in U$, there exists $n_x \in \Z_{\geq 0}$ such that 
$x+P_{n_x} \subset U$. 
We can directly check that $P$ is a topological group 
with respect to the additive structure. 
The following lemma gives a criterion 
for the continuity of the multiplication map. 

\begin{lem}\label{l-top-criterion}
Let $P$ be a pseudo-ring. 
Let $\{P_n\}_{n \in \Z_{\geq 0}}$ be a set of a pseudo-subrings of $P$ 
such that $P_0 \supset P_1 \supset \cdots$. 
We equip $P$ with the group topology induced by $\{P_n\}_{n \in \Z_{\geq 0}}$. 
Then the following are equivalent. 
\begin{enumerate}
\item 
$P$ is a topological pseudo-ring, i.e. the multiplication map 
$$\mu:P \times P \to P,\quad (x, y) \mapsto xy$$
is continuous. 
\item 
For any element $x \in P$ and any non-negative integer $n$, 
there exists a non-negative integer $m$ satisfying the inclusion  
$$
xP_m:=\{xy \in P\,|\, y \in P_m\} \subset P_n.
$$
\end{enumerate}
\end{lem}

\begin{proof}
Assume (1). 
Take an element $x \in P$ and a non-negative integer $n$. 
By (1), the composite map 
$$\theta:P \to P \times P \to P, \quad y \mapsto (x, y) \mapsto xy$$
is continuous. 
In particular, we get an inclusion $P_m \subset \theta^{-1}(P_n)$ 
for some $m \in \Z_{\geq 0}$. 
Therefore we obtain $x P_m=\theta(P_m) \subset \theta(\theta^{-1}(P_n)) \subset P_n$, 
hence (2) holds. 
Thus (1) implies (2).

Assume (2). 
Fix $z \in P$ and $n \in \Z_{> 0}$. 
It suffices to show that $\mu^{-1}(z+P_n)$ is 
an open subset of $P \times P$. 
If $\mu^{-1}(z+P_n)$ is empty, then there is nothing to show. 
Pick $(x, y) \in \mu^{-1}(z+P_n)$, i.e. $xy \in z+P_n$. 
By (2), we can find a positive integer $m$ 
such that $m \geq n$, $xP_m \subset P_n$ and $yP_m \subset P_n$. 
Take $x', y' \in P_m$. 
We get 
$$\mu(x+x', y+y')=xy+x'y+xy'+x'y' \in (z+P_n)+P_n+P_n+P_m \subset z+P_n.$$
Therefore it holds that 
$$(x+P_m) \times (y+P_m) \subset \mu^{-1}(z+P_n),$$ 
hence (1) holds. Thus (2) implies (1).  
\end{proof}

\subsection{Quotients by ideals}

The materials treated in this subsection 
appear in literature (cf. \cite[(1.4.1)]{Hub96}). 

\subsubsection{Quotients of f-adic rings}\label{ss-quot-fad}

Let $A$ be an f-adic ring. 
Let $J$ be an ideal of $A$ and let $\pi:A \to A/J$ be the natural ring homomorphism. 
For a ring of definition $A_0$ of $A$ and 
an ideal of definition $I_0$ of $A_0$, 
we equip $A/J$ with the group topology induced by $\{\pi(I_0^nA_0)\}_{n \in \Z_{>0}}$. 
We can check directly by definition that 
the topology of $A/J$ does not depend on the choices of $A_0$ and $I_0$. 
It follows from Lemma~\ref{l-top-criterion} that $A/J$ is a topological ring. 
Again by definition, 
we have that $A/J$ is an f-adic ring and $\pi:A \to A/J$ is adic. 
We call $A/J$ the {\em quotient f-adic ring} of $A$ by $J$. 

The following lemma should be well-known for experts, 
however we include the proof for the sake of completeness. 

\begin{lem}\label{l-2quot-top}
With the notation as above, the topology on $A/J$ coincides with the quotient topology induced by $\pi$.  
\end{lem}

\begin{proof}
To avoid confusion, we call the topology on $A/J$ defined above 
{\em the f-adic topology} in this proof. 
Fix a ring of definition $A_0$ of $A$ and an ideal of definition $I_0$ of $A_0$. 

Let $V$ be an open subset of $A/J$ with respect to the quotient topology. 
Pick $x \in \pi^{-1}(V)$. 
Since $\pi^{-1}(V)$ is an open subset of $A$, 
we can find a positive integer $n_x$ 
satisfying $x+I_0^{n_x}A_0 \subset \pi^{-1}(V)$. 
Thus we get an equation 
$\pi^{-1}(V)=\bigcup_{x \in \pi^{-1}(V)} (x+I_0^{n_x}A_0),$ which implies 
$$V=\pi\left(\pi^{-1}(V)\right)=\bigcup_{x \in \pi^{-1}(V)} \left(\pi(x)+\pi(I_0^{n_x}A_0)\right).$$
Therefore, $V$ is an open subset of $A/J$ with respect to the f-adic topology. 

Let $V$ be an open subset of $A/J$ with respect to the f-adic topology. 
For any $y \in V$, we can find a positive integer $n_y$ such that 
$y+\pi(I_0^{n_y}A_0) \subset V$. 
Thus we get an equation $V=\bigcup_{y \in V} (y+\pi(I_0^{n_y}A_0)),$ 
which implies 
$$\pi^{-1}(V)=\bigcup_{y \in V} \pi^{-1}(y+\pi(I_0^{n_y}A_0))
=\bigcup_{y \in V} (y'+J+I_0^{n_y}A_0),$$
where $y'$ is an element of $A$ satisfying $\pi(y')=y$. 
Thus $\pi^{-1}(V)$ is an open subset of $A$. 
Therefore, we have that $V$ is an open subset of $A/J$ with respect to the quotient topology. 
\end{proof}

\subsubsection{Quotients of affinoid rings}\label{ss-quot-aff}

Let $A=(A^{\rhd}, A^+)$ be an affinoid ring. 
For an ideal $J^{\rhd}$ of $A^{\rhd}$, 
we equip $A^{\rhd}/J^{\rhd}$ with the quotient topology. 
Thanks to Subsection~\ref{ss-quot-fad}, 
we have that $A^{\rhd}/J^{\rhd}$ is an f-adic ring and 
the induced ring homomorphism $\pi:A^{\rhd} \to A^{\rhd}/J^{\rhd}$ is adic. 
We set $(A/J^{\rhd})^+$ to be the integral closure of $\pi(A^+)$ in 
$A^{\rhd}/J^{\rhd}$. 
Then the pair $A/J^{\rhd}:=(A^{\rhd}/J^{\rhd}, (A/J^{\rhd})^+)$ 
is an affinoid ring and 
$\pi:A \to A/J^{\rhd}$ is an adic ring homomorphism of affinoid rings. 
We call $A/J^{\rhd}$ the {\em quotient affinoid ring} of $A$ by $J^{\rhd}$.

\begin{rem}
We use the same notation as above. 
Then $A/J^{\rhd}$ satisfies the following universal property: 
for any continuous ring homomorphism $\varphi:A \to B$ 
to an affinoid ring $B$ such that $\varphi(J^{\rhd})=0$, 
there exists a unique continuous ring homomorphism 
$\psi:A/J^{\rhd} \to B$ such that $\varphi=\psi \circ \pi$. 
\end{rem}

\subsection{Hausdorff quotient}

Let $A$ be a topological ring. 
Then the closure $\overline{\{0\}}$ of the zero ideal in $A$ 
is an ideal of $A$ and 
the residue ring $A^{\hd}:=A/\overline{\{0\}}$, 
equipped with the quotient topology, is a Hausdorff topological ring. 
Both of $A^{\hd}$ and $A \to A^{\hd}$ are 
called the {\em Hausdorff quotient} of $A$. 
For a continuous ring homomorphism $\varphi:A \to B$, 
we get a commutative diagram of continuous ring homomorphisms: 
$$\begin{CD}
A @>>> B\\
@VVV @VVV\\
A^{\hd} @>>> B^{\hd}.
\end{CD}$$

\begin{rem}\label{r-hd-univ}
Let $A$ be a topological ring and let $\varphi:A \to B$ 
be a continuous ring homomorphism to a Hausdorff topological ring $B$. 
Then there exists a unique continuous ring homomorphism 
$\psi:A^{\hd} \to B$ such that $\varphi=\psi \circ \pi$, 
where $\pi:A \to A^{\hd}$ is the natural map. 
\end{rem}

\begin{rem}\label{r-quot-open}
Let $A$ be a topological ring and 
let $\pi:A \to A^{\hd}$ be its Hausdorff quotient. 
Clearly, the induced map $\pi^*:\mathfrak O_{A^{\hd}} \to \mathfrak O_{A}$ is bijective and its inverse map $\pi_*:=(\pi^*)^{-1}$ 
satisfies $\pi_*(A_0)=\pi(A_0)$ for any $A_0 \in \mathfrak O_{A}$. 
Moreover we have that $\pi_*(A_0)$ is 
canonically isomorphic to $A_0^{\hd}$. 
\end{rem}

Let $A$ be an f-adic ring. 
For a ring of definition $A_0$ of $A$ and 
an ideal of definition $I_0$ of $A_0$, 
we have that $\overline{\{0\}} \subset I_0^nA_0$. 
Thus its image $\pi(I_0^nA_0)$ is an open subset of $A^{\hd}$. 
Since $\{\pi(I_0^nA_0)\}_{n \in \Z_{>0}}$ is a fundamental system 
of open neighbourhoods of $0 \in A^{\hd}$ by Lemma~\ref{l-2quot-top}, 
we have that $A^{\hd}$ is an f-adic ring. 
Moreover, the natural ring homomorphism $A \to A^{\hd}$ is adic.

\begin{prop}\label{p-hd-open}
Let $A$ be an f-adic ring and let $\pi:A \to A^{\hd}$ 
be the Hausdorff quotient of $A$. 
Then the following hold. 
\begin{enumerate}
\item 
It holds that $\pi_*(\mathfrak B_A)\subset \mathfrak B_{A^{\hd}}$ 
and 
$\pi^*(\mathfrak B_{A^{\hd}}) \subset \mathfrak B_A$. 
\item 
It holds that 
$\pi(A^{\circ})=(A^{\hd})^{\circ}$ and 
$\pi^{-1}((A^{\hd})^{\circ})=A^{\circ}$. 
\item 
It holds that $\pi_*(\mathfrak I_A)\subset \mathfrak I_{A^{\hd}}$ 
and 
$\pi^*(\mathfrak I_{A^{\hd}}) \subset \mathfrak I_A$.
\end{enumerate}
\end{prop}

\begin{proof}
We now show (1). 
The former inclusion $\pi_*(\mathfrak B_A)\subset \mathfrak B_{A^{\hd}}$ 
follows from \cite[Lemma 1.8(i)]{Hub93}. 
Let us prove the latter one $\pi^*(\mathfrak B_{A^{\hd}}) \subset \mathfrak B_A$. 
Take $B_1 \in \mathfrak B_{A^{\hd}}$. 
For any positive integer $n_1$, 
the boundedness of $B_1$ enables us 
to find a positive integer $n_2$ such that 
$$\pi(I_0^{n_2}A_0) \cdot B_1 \subset \pi(I_0^{n_1}A_0).$$
This inclusion immediately induces the following inclusion
$$I_0^{n_2}A_0 \cdot \pi^{-1}(B_1) \subset I_0^{n_1}A_0,$$
which implies $\pi^*(\mathfrak B_{A^{\hd}}) \subset \mathfrak B_A$. 
Thus (1) holds. 
The assertion (2) follows from (1) and \cite[Corollary 1.3(iii)]{Hub93}. 

We now show (3). 
Let us prove the former inclusion 
$\pi_*(\mathfrak I_A)\subset \mathfrak I_{A^{\hd}}$. 
Take $A^+ \in \mathfrak I_A$. 
It follows from (2) that $\pi(A^+) \subset (A^{\hd})^{\circ}$. 
Take an element $x \in A$ whose image 
$\pi(x)$ is integral over $\pi(A^+)$. 
We get 
$$x^n+a_1x^{n-1}+\cdots+a_n \in \overline{\{0\}}$$
for some $a_1, \cdots, a_n \in A^+$. 
By $\overline{\{0\}} \subset A^+$, we have that $x$ is integral over $A^+$. 
Since $A^+$ is integrally closed in $A$, 
we obtain $x \in A^+$. 
Hence, it holds that $\pi(x) \in \pi(A^+)$. 
Therefore, we get $\pi(A^+) \in \mathfrak I_{A^{\hd}}$ and 
$\pi_*(\mathfrak I_A) \subset \mathfrak I_{A^{\hd}}$. 
This completes the proof of the inclusion $\pi_*(\mathfrak I_A)\subset \mathfrak I_{A^{\hd}}$. 

Let us prove the other inclusion 
$\pi^*(\mathfrak I_{A^{\hd}}) \subset \mathfrak I_A$. 
Take $B^+ \in \mathfrak I_{A^{\hd}}$. 
It follows from (2) that 
$$\pi^{-1}(B^+) \subset \pi^{-1}((A^{\hd})^{\circ})=A^{\circ}.$$ 
Take an element $x \in A$ which is integral over $\pi^{-1}(B^+)$. 
We get 
$$x^n+a_1x^{n-1}+\cdots+a_n=0$$
for some $a_1, \cdots, a_n \in \pi^{-1}(B^+)$. 
Since $B^+$ is integrally closed in $A^{\hd}$, 
we have that $\pi(x) \in B^+$ and $x \in \pi^{-1}(B^+)$. 
Therefore, we get $\pi^{-1}(B^+) \in \mathfrak I_{A}$ and 
$\pi^*(\mathfrak I_{A^{\hd}}) \subset \mathfrak I_{A}$. 
Thus (3) holds. 
\end{proof}

\subsection{Completion}

Let $A$ be a topological ring. 
Its completion $A \to \widehat{A}$ 
is defined in \cite[Ch II \S 3 Section 3, Ch. III \S 6 Sect. 3]{Bou89} 
(cf. \cite[Ch. 0, Section 7.1(c)]{FK}). 
When we treat f-adic rings, 
their completions can be constructed 
by the classical method using Cauchy sequences 
(cf. \cite[Section 10.1]{AM69}). 
Even for general topological rings, 
we can construct the completions by using Cauchy nets in a similar way. 

If the topology of $A$ coincides with the $I$-adic topology for some finitely generated ideal $I$, 
then its completion $\widehat{A}$ is constructed 
also by the inverse limit $\varprojlim_n A/I^n$. 

\begin{rem}
Let $A$ be a ring, $I$ a finitely generated ideal of $A$, and $M$ an $A$-module. 
Then the inverse limit $\widehat{M}:=\varprojlim_n M/I^nM$ is complete 
with respect to the projective limit topology and $M/I^nM \xrightarrow{\simeq} \widehat{M}/I^n\widehat{M}$ 
by \cite[Ch. 0, Lemma 7.2.8 and Proposition 7.2.16]{FK}. 
If $I$ is not finitely generated, 
there exist counterexamples to these assertions (\cite[Ch. 0, Example 7.2.10]{FK}). 
\end{rem}

Let $A$ be an f-adic ring and let $\gamma:A \to \widehat{A}$ be the completion. 
For any open subring $A_0$ of $A$, 
the induced continuous ring homomorphsim $\widehat{A_0} \to \widehat{A}$ 
is injective and open, hence we consider $\widehat{A_0}$ 
as an open subring of $\widehat{A}$. 
In other words, we get a map 
$\gamma_1:\mathfrak O_A \to \mathfrak O_{\widehat{A}}$ 
defined by $\gamma_1(A_0)=\overline{\gamma(A_0)}$ for any $A_0 \in \mathfrak O_A$. 
For a ring of definition $A_0$ of $A$ and 
an ideal of definition $I_0$ of $A_0$, 
it follows from \cite[Lemma 1.6(i)]{Hub93} that $\{I_0^n\widehat{A_0}\}$ 
is a fundamental system of open neighbourhoods of $0 \in \widehat{A}$. 
In particular, 
the completion $\gamma:A \to \widehat{A}$ is adic. 
The completion $\gamma$ uniquely factors through the Hausdorff quotient: 
$$\gamma:A \xrightarrow{\pi}  A^{\hd} \xrightarrow{\gamma'}\widehat{A}.$$
Thanks to \cite[Corollary 1.9(ii)]{Hub93}, 
also the induced map $\gamma'$ is an adic ring homomorphism.

\begin{lem}\label{l-complete-open}
Let $A$ be an f-adic ring and 
let $\gamma:A \to \widehat{A}$ be the completion. 
Consider the following two maps: 
\begin{eqnarray*}
\gamma_1:\mathfrak O_{A} \to \mathfrak O_{\widehat{A}},&& 
A_0 \mapsto \overline{\gamma(A_0)}\\
\gamma^*:\mathfrak O_{\widehat{A}} \to \mathfrak O_{A},&& 
B_0 \mapsto \gamma^{-1}(B_0).
\end{eqnarray*}
Then both $\gamma_1 \circ \gamma^*$ and $\gamma^*\circ \gamma_1$ are 
the identity maps. 
\end{lem}

According to Definition~\ref{d-OBI}(1), we set $\gamma_*:=\gamma_1$. 

\begin{proof}
By Remark~\ref{r-quot-open}, 
we may assume that $A$ is Hausdorff and 
$A$ is a subring of $\widehat{A}$. 

First we show that $\gamma^* \circ \gamma_1$ is the identity map. 
Take $A_0 \in \mathfrak O_{A}$. 
We show $A_0=\overline{A_0} \cap A$. 
It suffices to prove that $A_0 \supset \overline{A_0} \cap A$. 
Fix $\alpha \in \overline{A_0} \cap A$. 
There is a Cauchy sequence 
$\{a_n\}_{n \in \Z_{>0}} \subset A_0$ converging to $\alpha \in A$. 
Since $A_0$ is a closed subset of $A$, it follows that $\alpha \in A_0$. 
Therefore we get $A_0 \supset \overline{A_0} \cap A$, 
hence $\gamma^* \circ \gamma_1$ is the identity map.

Second we show that $\gamma_1 \circ \gamma^*$ is the identity map. 
Take $B_0 \in \mathfrak O_{\widehat{A}}$. 
We prove $B_0=\overline{B_0 \cap A}$. 
It suffices to show that $B_0 \subset \overline{B_0 \cap A}$. 
Take $\beta \in B_0$. 
There exists a Cauchy sequence 
$\{b_n\}_{n \in \Z_{>0}} \subset A$ converging to $\beta$. 
Since $B_0$ is an open subset of $\widehat{A}$, we can find $N \in \Z_{>0}$ 
such that $b_n-\beta \in B_0$ for any $n \geq  N$. 
Thus we get $b_n=\beta+(b_n-\beta) \in B_0$ for any $n \in \Z_{\geq N}$. 
In particular, the shifted sequence $b_N, b_{N+1}, \cdots$ 
is contained in $A \cap B_0$. 
Therefore, we get $\beta \in \overline{B_0 \cap A}$, 
hence $\gamma_1 \circ \gamma^*$ is the identity map. 
\end{proof}

\begin{lem}\label{l-sub-bdd-bdd}
Let $B$ be an f-adic ring and 
let $A$ be an f-adic subring of $B$ 
whose topology coincides with the induced topology from $B$. 
Set $j:A \hookrightarrow B$ to be the inclusion map. 
Then the inclusion 
$j^*(\mathfrak B_{B}) \subset \mathfrak B_{A}$ holds. 
\end{lem}

\begin{proof}
Let $B_0$ be a bounded subset of $B$. 
It suffices to show that $B_0 \cap A$ is a bounded subset of $A$. 
Let $U$ be an open neighbourhood of $0 \in A$. 
There exists an open subset $U_B$ of $B$ such that $U_B \cap A=U$. 
Since $B_0$ is bounded, 
we can find an open neighbourhood $V_B$ of $0 \in B$ such that 
$B_0 \cdot V_B \subset U_B$. 
We get  
$$(B_0 \cap A)\cdot (V_B \cap A) \subset 
(B_0 \cdot V_B) \cap A \subset U_B \cap A=U,$$
which implies that $B_0 \cap A$ is bounded. 
\end{proof}

\begin{lem}\label{l-complete-open2}
Let $A$ be an f-adic ring and 
let $\gamma:A \to \widehat{A}$ be the completion of $A$. 
Let $\gamma_*:\mathfrak O_A \to \mathfrak O_{\widehat{A}}$ and $\gamma^*:\mathfrak O_{\widehat{A}} \to \mathfrak O_{A}$ be 
the induced maps (cf. Lemma~\ref{l-complete-open}). 
Then the following assertions hold. 
\begin{enumerate}
\item If $A_0 \in \mathfrak B_{A}$, then $\gamma_*(A_0) 
\in \mathfrak B_{\widehat{A}}$.
\item If $B_0 \in \mathfrak B_{\widehat{A}}$, then $\gamma^*(B_0) \in \mathfrak B_{A}$. 
\end{enumerate}
\end{lem}

\begin{proof}
By Proposition~\ref{p-hd-open}, we may assume that $A$ is Hausdorff and 
$A$ is a subring of $\widehat{A}$. 
The assertion (1) follows from the fact that 
$\gamma_*(A_0)=\overline{A_0}$ is bounded (\cite[Lemma 1.6(i)]{Hub93}). 
The assertion (2) holds by Lemma~\ref{l-sub-bdd-bdd}. 
\end{proof}

\begin{lem}\label{l-pbdd-completion}
Let $A$ be an f-adic ring. 
Then the natural map $\widehat{A^{\circ}} \to (\widehat{A})^{\circ}$ is an isomorphism of topological rings.
\end{lem}

\begin{proof}
Let $\gamma:A \to \widehat{A}$ be the completion of $A$. 

Note that the natural map $\varphi:\widehat{A^{\circ}} \to (\widehat{A})^{\circ}$ in the statement is 
constructed as follows. 
Since the completion $\gamma:A \to \widehat{A}$ is adic by 
\cite[Lemma 1.6]{Hub93}, 
it follows from \cite[Corollary 1.3(iii), Lemma 1.8(i)]{Hub93} that $\gamma(A^{\circ}) \subset (\widehat{A})^{\circ}$. 
Taking the completion of $\gamma|_{A^{\circ}}:A^{\circ} \to (\widehat{A})^{\circ}$, 
we obtain a natural continuous ring homomorphism 
$\varphi:\widehat{A^{\circ}} \to \widehat{(\widehat{A})^{\circ}}=(\widehat{A})^{\circ}$, 
where the equation $\widehat{(\widehat{A})^{\circ}}=(\widehat{A})^{\circ}$ 
follows from the fact that $(\widehat{A})^{\circ}$ is 
an open, hence a closed, subring of $\widehat{A}$.

By the construction of $\varphi$, 
it follows that $\varphi$ is open and injective, 
hence we can consider $\widehat{A^{\circ}}$ as an open subring of $(\widehat{A})^{\circ}$. 
It suffices to show $\widehat{A^{\circ}} \supset (\widehat{A})^{\circ}$. 
Take $\alpha \in (\widehat{A})^{\circ}$. 
Let $B_0$ be a bounded open subring of $(\widehat{A})^{\circ}$ containing $\alpha$. 
We set $A_0:=\gamma^{-1}(B_0)$. 
Then $A_0$ is a bounded open subring of $A$ such that  $\widehat{A_0}=B_0$ by Lemma~\ref{l-complete-open2}. 
Thus we obtain $\alpha \in \widehat{A_0} \subset \widehat{A^{\circ}}$, 
as desired. 
\end{proof}

\begin{lem}\label{l-integral-completion}
Let $B$ be an f-adic ring and let $A$ be an open subring of $B$. 
Let $\widehat{A}$ and $\widehat{B}$ be the completions of $A$ and $B$, 
respectively. 
Set $\gamma:B \to \widehat{B}$ be the induced map. 
Then the following hold. 
\begin{enumerate}
\item The equation $\widehat{B}=\gamma(B) \cdot \widehat{A}$ holds, 
where $\gamma(B) \cdot \widehat{A}$ denotes the smallest subring of $\widehat{B}$ containing $\gamma(B) \cup \widehat{A}$. 
\item If $A \subset B$ is an integral extension, then so is $\widehat{A} \subset \widehat{B}$. 
\end{enumerate}
\end{lem}

\begin{proof}
We first show (1). Set $C:=\gamma(B) \cdot \widehat{A}$. 
Since $\widehat{A}$ is an open subring of $\widehat{B}$, 
so is $C$. 
In particular, $C$ is a closed subset of $\widehat{B}$ containing $\gamma(B)$. 
Thus $\widehat{B}=C$, hence (1) holds. 
The assertion (2) immediately follows from (1). 
\end{proof}

\begin{lem}\label{l-complete-open3}
Let $A$ be an f-adic ring. 
Let $\gamma:A \to \widehat{A}$ be the completion and 
let $\gamma_*:\mathfrak O_A \to \mathfrak O_{\widehat{A}}$ and $\gamma^*:\mathfrak O_{\widehat{A}} \to \mathfrak O_{A}$ be 
the induced maps (cf. Lemma~\ref{l-complete-open}). 
Then the following hold. 
\begin{enumerate}
\item If $A^+ \in \mathfrak I_{A}$, then $\gamma_*(A^+) 
\in \mathfrak I_{\widehat{A}}$.
\item If $B^+ \in \mathfrak I_{\widehat{A}}$, then $\gamma^*(B^+) \in \mathfrak I_{A}$. 
\end{enumerate}
\end{lem}

\begin{proof} 
By Proposition~\ref{p-hd-open}, 
we may assume that $A$ is Hausdorff 
and $A$ is an subring of $\widehat{A}$. 

We first show (1). 
Since $\gamma_*(A^+)=\widehat{A^+}$ and 
$\widehat{A^+} \subset \widehat{A^{\circ}} =(\widehat{A})^{\circ}$ 
(Lemma~\ref{l-pbdd-completion}), 
it suffices to show that $\widehat{A^+}$ 
is integrally closed in  $\widehat{A}$. 
Let $\beta \in \widehat{A}$ be an element which is integral over $\widehat{A^+}$. 
We can write 
$$\beta^n+\alpha_1\beta^{n-1}+\cdots+\alpha_n=0$$
for some $\alpha_i \in \widehat{A^+}$. 
There exist sequences $\{a_{i, k}\}_{k \in \Z_{>0}}\subset A^+$ and 
$\{b_k\}_{k \in \Z_{>0}} \subset A$ converging to $\alpha_i$ and $\beta$, 
respectively. 
We set 
$$b_k^n+a_{1, k}b_k^{n-1}+\cdots+a_{n, k}=:c_k \in A.$$
We see that the sequence $c_1, c_2, \cdots$ converges to zero. 
Since $A^+$ is an open subset of $A$, 
we can assume that $\{c_k\}_{k \in \Z_{>0}} \subset A^+$. 
In particular, each $b_k$ is integral over $A^+$. 
As $A^+$ is integrally closed in $A$, 
it follows that $\{b_k\}_{k \in \Z_{>0}} \subset A^+$. 
Therefore, we get $\beta \in \widehat{A^+}$, hence (1) holds.

We now show (2). 
Set $\gamma^*(B^+)=\gamma^{-1}(B^+)=:A_1$. 
Since Lemma~\ref{l-pbdd-completion} implies 
$$B^+ \subset (\widehat{A})^{\circ}=\widehat{A^{\circ}}=\gamma_*(A^{\circ}),$$
we have that $A_1=\gamma^*(B^+) \subset \gamma^*\gamma_*(A^{\circ})=A^{\circ}$. 
Thus it suffices to prove that $A_1$ is integrally closed in $A$. 
Let $A_2$ be the integral closure of $A_1$ in $A$. 
It follows from $A_1 \subset A_2$ that $A_2$ is an open subring of $A$. 
Since $A_2$ is integral over $A_1$, 
it holds by Lemma \ref{l-integral-completion}(2) that 
$\gamma_*(A_2)=\overline{A_2}$ is integral over $\gamma_*(A_1)=\gamma_*\gamma^*(B^+)=B^+$. 
Thanks to $B^+ \in \mathfrak I_{\widehat{A}}$, 
we get $\overline{A_2} \subset B^+$. 
By Lemma~\ref{l-complete-open}, we have that  
$$A_2 =\gamma^*\gamma_*A_2=\gamma^{-1}(\overline{A_2}) \subset \gamma^{-1}(B^+)=A_1,$$ 
as desired. Thus (2) holds. 
\end{proof}

\begin{lem}\label{l-common-comp}
Let $A$ be a topological ring and let $\theta:A \to \widehat A$ be the completion of $A$. 
Let $\varphi:A \to B$ be a continuous ring homomorphism 
to a topological ring $B$ that satisfies the following properties:  
\begin{enumerate}
\item There exists a continuous 
ring homomorphism $\psi:B \to \widehat{A}$ 
such that $\theta=\psi \circ \varphi$. 
\item For any open neighbourhood $V$ of $0 \in B$, 
there exists an open neighbourhood $W$ of $0 \in \widehat{A}$ such that $\psi^{-1}(W) \subset V$. 
\end{enumerate}
Then the induced map $\widehat{\varphi}:\widehat{A} \to \widehat{B}$ 
is an isomorphism of topological rings. 
\end{lem}

\begin{proof}
We set $\theta_A:=\theta$ and let 
$\theta_B:B \to \widehat{B}$ be the completion, 
so that we get a commutative diagram: 
$$\begin{CD}
A @>\varphi >> B @>\psi >> \widehat{A}\\
@VV\theta_A V @VV\theta_B V @V\simeq V \widehat{\theta_A}V\\
\widehat{A} @>\widehat{\varphi}>> \widehat{B} @>\widehat{\psi} >> \widehat{\widehat{A}}.
\end{CD}$$

\begin{step}\label{s-common-comp1}
$\widehat{\psi}$ is injective. 
\end{step}

\begin{proof}(of Step~\ref{s-common-comp1}) 
Take $\beta \in \widehat{B}$ such that $\widehat{\psi}(\beta)=0$. 
There is a Cauchy net $\{b_i\}_{i \in I} \subset B$ whose image 
$\{\theta_B(b_i)\}_{i \in I}$ converges to $\beta$. 
Since $\{\psi(b_i)\}_{i \in I}$ is a Cauchy net of $\widehat{A}$, 
it converges to an element $\alpha$ of $\widehat{A}$. 
We have that 
$$0=\widehat{\psi}(\beta)=\widehat{\psi}(\lim_{i \in I} \theta_{B}(b_i))
=\widehat{\theta_A}(\alpha),$$
which implies $\alpha=0$. 
Take an open neighbourhood $V$ of $0 \in B$. 
It follows from (2) that there exists an open neighbourhood $W$ 
of $0 \in \widehat{A}$ such that $\psi^{-1}(W) \subset V$. 
Since $0=\alpha=\varinjlim_{i \in I}\psi(b_i)$, 
there exists an index $i_0 \in I$ such that $\psi(b_i) \in W$ 
for any $i \in I_{\geq i_0}$. 
In particular, 
we get $b_i \in \psi^{-1}(W) \subset V$ for any $i \in I_{\geq i_0}$. 
This implies that $\{b_i\}_{i \in I}$ converges to zero, 
hence $\beta=\theta_B(\lim_{i \in I} b_i)=0$, as desired. 
This completes the proof of Step~\ref{s-common-comp1}. 
\end{proof}

\begin{step}\label{s-common-comp2}
$\theta_B=\widehat{\varphi} \circ \psi$.
\end{step}

\begin{proof}(of Step~\ref{s-common-comp2}) 
Since $\widehat{\psi}$ is injective by Step~\ref{s-common-comp1}, 
it suffices to show that 
$\widehat{\psi} \circ \theta_B=\widehat{\psi} \circ 
\widehat{\varphi} \circ \psi$, which follows from 
$$\widehat{\psi} \circ 
\widehat{\varphi} \circ \psi=\widehat{\theta_A} \circ \psi=\widehat{\psi} \circ \theta_B,$$
where the first equation is guaranteed by (1). 
This completes the proof of Step~\ref{s-common-comp2}.
\end{proof}

\begin{step}\label{s-common-comp3}
$\widehat{\varphi}$ is bijective. 
\end{step}

\begin{proof}(of Step~\ref{s-common-comp3}) 
Since $\widehat{\psi} \circ \widehat{\varphi}$ is bijective, 
we have that $\widehat{\varphi}$ is injective. 

It suffices to show that $\widehat{\varphi}$ is surjective. 
Take $\beta \in \widehat{B}$. 
We can find a Cauchy net $\{b_i\}_{i \in I} \subset B$ 
such that its image $\{\theta_B(b_i)\}_{i \in I}$ converges to $\beta$. 
Since $\psi:B \to \widehat{A}$ is continuous, 
also $\{\psi(b_i)\}_{i \in I}$ is a Cauchy net. 
Thus we can find an element $\alpha \in \widehat{A}$ 
with $\alpha=\lim_{i \in I}\psi(b_i)$. 
We have that 
$$\widehat{\varphi}(\alpha)=\widehat{\varphi}\left(\lim_{i \in I}\psi(b_i)\right)=\lim_{i \in I}\theta_B(b_i)=\beta,$$
where the second equation holds by Step~\ref{s-common-comp2}. 
Thus $\widehat{\varphi}$ is surjective. 
This completes the proof of Step~\ref{s-common-comp3}.
\end{proof}

\begin{step}\label{s-common-comp4}
$\widehat{\varphi}$ is an open map. 
\end{step}

\begin{proof}(of Step~\ref{s-common-comp4}) 
Let $U$ be an open subset of $\widehat{A}$. 
Since both $\widehat{\varphi}$ and $\widehat{\psi}$ are bijective by Step~\ref{s-common-comp3}, 
we get an equation: 
$$\widehat{\varphi}(U)=\widehat{\psi}^{-1}(\widehat{\psi}(\widehat{\varphi}(U))).$$
Since $\widehat{\psi}(\widehat{\varphi}(U))$ is an open subset and 
$\widehat{\psi}$ is continuous, 
$\widehat{\varphi}(U)$ is an open set. 
This completes the proof of Step~\ref{s-common-comp4}.
\end{proof}
Step~\ref{s-common-comp3} 
and Step~\ref{s-common-comp4} 
complete the proof of Lemma~\ref{l-common-comp}. 
\end{proof}

\subsection{Topological tensor products}

In this subsection, we introduce tensor products 
for the category of f-adic rings (Theorem~\ref{t-top-tensor}).

\begin{dfn}\label{d-tensor}
Let $\mathcal R$ be a category. 
For two arrows $\varphi:R \to A$ and $\psi:R \to B$ in $\mathcal R$, 
we say that $(T, f, g)$ is a {\em tensor product} in $\mathcal R$ of $(\varphi, \psi)$ or of a diagram 
$$\begin{CD}
R @>\varphi >> A\\
@VV\psi V\\
B,
\end{CD}$$
if the following hold: 
\begin{enumerate}
\item 
$T$ is an object of $\mathcal R$, 
\item 
$f:A \to T$ and $g:B \to T$ are arrows of $\mathcal R$ 
such that $f \circ \varphi=g \circ \psi$, and
\item 
for a commutative diagram of arrows of $\mathcal R$ 
$$\begin{CD}
R @>\varphi >> A\\
@VV\psi V @VVf'V\\
B @>g'>>C,
\end{CD}$$
there exists a unique arrow $\theta:T \to C$ in $\mathcal R$ 
satisfying $f'=\theta \circ f$ and $g'=\theta \circ g'$. 
\end{enumerate}
We often call $T$ {\em the tensor product} of $(\varphi, \psi)$ 
in $\mathcal R$ if no confusion arises. 
\end{dfn}

\begin{dfn}\label{d-ring-tensor}
The tensor products in the category of rings are called 
{\em ring-theoretic tensor products}. 
\end{dfn}

\begin{dfn}\label{d-categories}
\begin{enumerate}
\item 
Let $(\text{TopRing})$ be the category of topological rings  
whose arrows are the continuous ring homomorphisms. 
\item 
Let $(\text{FadRing})$ be the full subcategory of $(\text{TopRing})$ 
whose objects are f-adic rings.  
\item 
Let $(\text{FadRing})^{\text{ad}}$ 
be the category of f-adic rings 
whose arrows are the adic ring homomorphisms. 
\end{enumerate}
\end{dfn}

\begin{lem}\label{l-underlying}
Let $\mathcal R$ be one of 
the categories 
$({\rm TopRing})$ and $({\rm FadRing})$. 
Let $\varphi:R \to A$ and $\psi:R \to B$ be two arrows 
of $\mathcal R$. 
Assume that there exists a tensor product $(T, f, g)$ in $\mathcal R$ 
of $(\varphi, \psi)$, 
then the induced ring homomorphism $A \otimes_R B \to T$ 
from the ring-theoretic tensor product $A \otimes_R B$ is bijective. 
\end{lem}

\begin{proof}
We only treat the case where $\mathcal R=(\text{FadRing})$, 
as both the proofs are the same. 
Take a ring $C$ and a commutative diagram of ring homomorphisms: 
$$\begin{CD}
R @>\varphi >> A\\
@VV\psi V @VVf'V\\
B @>g'>> C.
\end{CD}$$
We equip $C$ with the trivial topology. 
Then $C$ is an f-adic ring such that $f'$ and $g'$ are continuous. 
Since $(T, f, g)$ is a tensor product in $(\text{FadRing})$ of $(\varphi, \psi)$, 
there exists a unique continuous ring homomorphism 
$\theta:T \to C$ such that $f'=\theta \circ f$ and $g'=\theta \circ g$. 
Take another ring homomorphism $\widetilde{\theta}:T \to C$ 
satisfying $f'=\widetilde \theta \circ f$ and $g'=\widetilde \theta \circ g'$. 
Then $\widetilde \theta$ is automatically continuous, 
hence it follows from Definition~\ref{d-tensor} that 
$\theta=\widetilde \theta$. 
Therefore the underlying ring $T$ satisfies the universal property 
characterising ring-theoretic tensor products. 
Thus the induced ring homomorphism $A \otimes_R B \to T$ 
is bijective. 
\end{proof}

\begin{thm}\label{t-top-tensor}
Let $\varphi:R \to A$ and $\psi:R \to B$ 
be adic ring homomorphisms of f-adic rings. 
Then the following hold:  
\begin{enumerate}
\item 
There exists a tensor product $(T, f, g)$ in $({\rm FadRing})$ 
of $(\varphi, \psi)$. 
\item 
$(T, f, g)$ is a tensor product in $({\rm FadRing})^{{\rm ad}}$ 
of $(\varphi, \psi)$. 
\item 
The induced ring homomorphism $A \otimes_R B \to T$ 
from the ring-theoretic tensor product $A \otimes_R B$ is bijective. 
\item 
The induced continuous ring homomorphisms $f:A \to T$ and $g:B \to T$ are adic. 
\end{enumerate}
\end{thm}

An f-adic ring $T$ satisfying the above properties 
is called a {\em topological tensor product} of $(\varphi, \psi)$. 
Because of (3), we often denote it by $A \otimes_R B$.

\begin{proof}
As a ring, we set $T:=A \otimes_R B$. 
Let 
$$\begin{CD}
R @>\varphi >> A\\
@VV\psi V @VVfV\\
B @>g>> T
\end{CD}$$
be the induced commutative diagram of ring homomorphisms. 
Let $A_0$ and $B_0$ be rings of definition of $A$ and $B$, respectively. Since $\varphi^{-1}(A_0) \cap \psi^{-1}(B_0)$ is an open subring of $R$, 
we can find a ring of definition $R_0$ of $R$ 
such that $R_0 \subset \varphi^{-1}(A_0) \cap \psi^{-1}(B_0)$. 
We have the natural ring homomorphism 
$$\rho:A_0 \otimes_{R_0} B_0 \to A \otimes_R B.$$
We take an ideal of definition $I_0$ of $R_0$. 
It follows from \cite[Lemma 1.8(ii)]{Hub93} that $I_0A_0$ and $I_0B_0$ are 
ideals of definition of $A_0$ and $B_0$, respectively. 
We equip $T=A \otimes_R B$ with the group topology defined by 
$$\left\{\rho\left(I_0^k\cdot(A_0 \otimes_{R_0} B_0)\right)\right\}_{k \in \Z_{>0}}.$$ 
We can check that this topology does not depend 
on the choices of $R_0, A_0,  B_0$ and $I_0$.

\begin{claim}
The ring $T$, equipped with the group topology defined as above, 
is an f-adic ring such that the induced ring homomorphisms $f:A \to T$ and $g:B \to T$ 
are adic.  
\end{claim}

\begin{proof}(of Claim) 
To show that $T$ is a topological ring, 
we prove that the property (2) of Lemma~\ref{l-top-criterion} holds. 
Take $\xi \in A \otimes_R B$ and $n \in \Z_{>0}$. 
We can write $\xi=\sum_{i=1}^r a_i \otimes_R b_i$ 
for some $a_i \in A$ and $b_i \in B$. 
Since $A$ and $B$ are topological rings, 
Lemma~\ref{l-top-criterion} enables us to find $m_1 \in \Z_{>0}$ satisfying  
$$a_i I_0^{m_1}A_0 \subset I_0^nA_0,\quad b_i I_0^{m_1}B_0 \subset I_0^nB_0$$ 
for any $i \in \{1, \cdots, r\}$. 
In particular, we obtain 
$$\xi \cdot \rho\left(I_0^{2m_1}\cdot(A_0 \otimes_{R_0} B_0)\right) 
\subset \rho\left(I_0^{2n}\cdot(A_0 \otimes_{R_0} B_0)\right),$$
hence (2) of Lemma~\ref{l-top-criterion} holds for $T$. 
Therefore, $T$ is a topological ring. 

Moreover, it follows from definition of the topology on $T$ 
that $T$ is an f-adic ring and both of $f$ and $g$ are adic. 
This completes the proof of Claim. 
\end{proof}

We will show the following property: 
\begin{enumerate}
\item[$(1)'$] $(T, f, g)$ is a tensor product in $({\rm FadRing})$ of $(\varphi, \psi)$. 
\end{enumerate}

For the time being, 
let us check that the statement of the theorem holds if $(1)'$ holds. 
We have that (1) follows from $(1)'$. 
We obtain (3) and (4) by the construction of $T$ and Claim, respectively. 
It follows from $(1)'$ and \cite[Corollary 1.9(ii)]{Hub93} that (2) automatically holds. 
Therefore, the statement of Theorem \ref{t-top-tensor} holds if $(1)'$ holds. 

\medskip 

Thus it suffices to show $(1)'$. 
Take a commutative diagram 
$$\begin{CD}
R @>\varphi >> A\\
@VV\psi V @VVf'V\\
B @>g'>> C 
\end{CD}$$
in $({\rm FadRing}).$ 
Since $T$ is a ring-theoretic tensor product, 
there exists a unique ring homomorphism $\theta:T \to C$ 
such that $f'=\theta \circ f$ and $g'=\theta \circ g$. 
It is enough to show that $\theta$ is continuous. 

Let $C_1$ be an open pseudo-subring of $C$. 
Since $f'$ and $g'$ are continuous, 
we can find a positive integer $k$ such that 
$I_0^kA_0 \subset f'^{-1}(C_1)$ and $I_0^kB_0 \subset g'^{-1}(C_1)$. 
Then we have that 
$$\{x \otimes_R 1\,|\, x \in I_0^kA_0\} \cup 
\{1 \otimes_R y\,|\, y \in I_0^kB_0\} \subset \theta^{-1}(C_1).$$ 
Since $\theta^{-1}(C_1)$ is a pseudo-subring of $T$, 
we get 
$$\rho(I_0^{2k}(A_0 \otimes_{R_0} B_0)) \subset \theta^{-1}(C_1).$$ 
Therefore $\theta$ is continuous. 
Thus $(1)'$ holds, 
hence so does the statement of Theorem~\ref{t-top-tensor}. 
\end{proof}

\subsection{Structure sheaves of Zariskian schemes}

\subsubsection{Affine case}\label{ss-FK-affine}

In this subsection, we summarise some results from \cite{FK} for later use. 
Let $A$ be a ring and let $I$ be an ideal of $A$. 
For an $A$-module $M$, recall that $\widetilde M$ 
is the quasi-coherent sheaf on $\Spec\,A$ satisfying 
$M=\Gamma(\Spec\,A, \widetilde M)$. 
We define a sheaf $M^{\diamondsuit}$ on a topological space $V(I)$ by 
$$M^{\diamondsuit}:=i^{-1}(\widetilde M).$$
The following theorem is nothing but \cite[Ch. I, Proposition B.1.4]{FK}, 
however we give a proof of it 
since the proof of \cite[Ch. I, Proposition B.1.4]{FK} omits some of arguments.

\begin{thm}\label{t-FK-sheafy}
Let $A$ be a ring and let $I$ be an ideal of $A$. 
For any $A$-module $M$ and element $f \in A$, 
the equation 
$$\Gamma(V(I) \cap D(f), M^{\diamondsuit})=M \otimes_A (1+IA_f)^{-1}A_f$$
holds. 
\end{thm}

\begin{proof}
Take elements $f_1, \cdots, f_n \in A$ such that 
$V(I) \subset D(f_1) \cup \cdots \cup D(f_n)$. 
For any $A$-module $N$ and any element $g \in A$, we set 
$$N_g^Z:=N \otimes_A (1+IA_g)^{-1}A_g$$
and $N^Z:=N_1^Z$. 
By the proof of \cite[Ch. I, Proposition B.1.4]{FK}, 
it suffices to show that the sequence 
\begin{equation}\label{zar-sch-sheafy}
0 \to M^Z \xrightarrow{\varphi} 
\prod_{1 \leq i\leq r} M_{f_i}^{Z} \xrightarrow{\psi} 
\prod_{1 \leq i<j \leq r} M_{f_if_j}^{Z}
\end{equation}
is exact, where $\varphi$ is defined by $\varphi(m)=(m \otimes_A 1, \cdots, m \otimes_A 1)$ for any $m \in M$ and 
$\psi$ is defined by the difference.

Replacing $A$ and $M$ by $A^Z$ and $M^Z$ respectively, 
we may assume that $A=(1+I)^{-1}A$. 
Thanks to the inclusion $V(I) \subset D(f_1) \cup \cdots \cup D(f_r)$, 
the images of the elements $f_1, \cdots, f_r$ to $A/I$ 
generate $A/I$. 
In particular, we can find elements $g_1, \cdots, g_r \in A$ such that 
$\sum_{i=1}^rg_if_i=1+z$ for some $z \in I$. 
Since $1+z \in 1+I \subset A^{\times}$, 
we may assume that $z=0$, i.e. the equation 
\begin{equation}\label{e-FK-sheafy1}
\sum_{i=1}^rg_if_i=1
\end{equation}
holds in $A$.  

\setcounter{step}{0}
\begin{step}\label{s-domain}
The sequence (\ref{zar-sch-sheafy}) is exact if $M=A/\p$ for some prime ideal of $A$. 
\end{step}

\begin{proof}(of Step~\ref{s-domain}) 
We may assume that $M=A$ and $A$ is an integral domain. 
Since $A$ is an integral domain, it is clear that $\varphi$ is injective. 
Take $(\xi_1, \cdots, \xi_r) \in \prod_{1 \leq i\leq r}(1+IA_{f_i})^{-1}A_{f_i}$ 
such that $\psi((\xi_1, \cdots, \xi_r))=0$. 
For each $i \in \{1, \cdots, r\}$, we can write 
$$\xi_i=\frac{\frac{a_i}{f_i^{n_i}}}{1+\frac{x_i}{f_i^{m_i}}}.$$
for some $a_i\in A$, $x_i \in I$ and $n_i, m_i \in \Z_{> 0}$. 
We may assume that there is a positive integer $n$ such that $n=n_i=m_i$ 
for any $i \in \{1, \cdots, r\}$. 
Moreover, replacing $f_i^n$ by $f_i$, 
the problem is reduced to the case where $n=1$. 
Thus we obtain 
$$\xi_i=\frac{\frac{a_i}{f_i}}{1+\frac{x_i}{f_i}}.$$
Since $\psi((\xi_1, \cdots, \xi_r))=0$, we get an equation
$$\frac{\frac{a_i}{f_i}}{1+\frac{x_i}{f_i}}=\frac{\frac{a_j}{f_j}}{1+\frac{x_j}{f_j}}$$
in $(1+IA_{f_if_j})^{-1}A_{f_if_j}$ for any $i, j\in \{1, \cdots, r\}$. 
Since $A$ is an integral domain, we get an equation 
\begin{equation}\label{e-FK-sheafy2}
(f_j+x_j)a_i=(f_i+x_i)a_j
\end{equation}
in $A$ for any $i, j\in \{1, \cdots, r\}$. 
We set 
$$a:=\sum_{j=1}^ra_jg_j.$$
Then, it holds that 
$$(f_i+x_i)a=\sum_{j=1}^r(f_i+x_i)a_jg_j=\sum_{j=1}^r(f_j+x_j)a_ig_j
=a_i\left(1+\sum_{j=1}^rx_jg_j\right),$$
where the second and third equations hold by 
(\ref{e-FK-sheafy2}) and (\ref{e-FK-sheafy1}) respectively. 
Since $1+\sum_{j=1}^rx_jg_j \in 1+IA \subset A^{\times}$, 
we get $(1+\sum_{j=1}^rx_jg_j)^{-1} \in A$. 
Thus, for 
$$\widetilde a:=\left(1+\sum_{j=1}^r x_jg_j\right)^{-1}a,$$
it holds that $\varphi(\widetilde a)=(\xi_1, \cdots, \xi_r)$, as desired. 
This completes the proof of Step~\ref{s-domain}. 
\end{proof}

\begin{step}\label{s-thickening}
Let $0 \to M_1 \to M_2 \to M_3 \to 0$ be an exact sequence of $A$-modules. 
If 
(\ref{zar-sch-sheafy}) is exact for the cases $M =M_1$ and $M=M_3$, then 
(\ref{zar-sch-sheafy}) is exact also for the case $M=M_2$. 
\end{step}

\begin{proof}(of Step~\ref{s-thickening}) 
We have the following commutative diagram 
$$\begin{CD}
@. 0 @. 0 @. 0 @.\\
@. @VVV @VVV @VVV\\
0 @>>> M_1^{Z} @>>> \prod_{1\leq i \leq r} (M_1)_{f_i}^Z @>>> \prod_{1 \leq i<j\leq r} (M_1)_{f_if_j}^Z\\
@. @VVV @VVV @VVV\\
0 @>>> M_2^{Z} @>>> \prod_{1\leq i \leq r} (M_2)_{f_i}^Z @>>> \prod_{1 \leq i<j\leq r} (M_2)_{f_if_j}^Z\\
@. @VVV @VVV @VVV\\
0 @>>> M_3^{Z} @>>> \prod_{1\leq i \leq r} (M_3)_{f_i}^Z @>>> \prod_{1 \leq i<j\leq r} (M_3)_{f_if_j}^Z\\
@. @VVV @VVV @VVV\\
@. 0 @. 0 @. 0 @.\\
\end{CD}$$
Since $A \to A^{Z}_f=(1+IA_f)^{-1}A_f$ is flat for any element $f \in A$, 
all the vertical sequences are exact. 
Moreover, also the upper and lower horizontal sequences are exact by assumption. 
It follows from a diagram chase that the middle horizontal sequence is exact, as desired. 
This completes the proof of Step~\ref{s-thickening}.
\end{proof}

\begin{step}\label{s-noether}
The sequence (\ref{zar-sch-sheafy}) is exact if $A$ is a noetherian ring and $M$ is a finitely generated $A$-module. 
\end{step}

\begin{proof}(of Step~\ref{s-noether}) 
Thanks to \cite[Theorem 6.5]{Mat89}, 
we can find a descending sequence of $A$-submodules of $M$: 
$$M=:M_0 \supset M_1 \supset \cdots \supset M_{\ell-1} \supset M_{\ell}=0$$
such that for any $k \in \{0, \cdots, \ell-1\}$, 
there exists a prime ideal $\p_k$ of $A$ such that $M_k/M_{k+1} \simeq A/\p_k$. 
By Step~\ref{s-thickening}, we may assume that $M=A/\p$ for some prime ideal $\p$ of $A$. 
The assertion of Step~\ref{s-noether} holds by Step~\ref{s-domain}. 
\end{proof}

\begin{step}\label{s-general}
The sequence (\ref{zar-sch-sheafy}) is exact without any additional assumptions. 
\end{step}

\begin{proof}(of Step~\ref{s-general}) 
We only show the exactness on $\prod_{1 \leq i\leq r} M_{f_i}^{Z}$, as the remaining case is easier. 
Take $(\xi_1, \cdots, \xi_r) \in \prod_{1 \leq i\leq r} M_{f_i}^{Z}$ 
such that $\psi((\xi_1, \cdots, \xi_r))=0$. 
Since 
$$M^Z_{f_i}=M \otimes_A (1+IA_{f_i})^{-1}A_{f_i},$$
we can write 
$$\xi_i=\frac{\frac{\mu_i}{f_i^{n_i}}}{1+\frac{x_i}{f_i^{m_i}}}$$
for some $x_i \in I$, $\mu_i \in M$, $n_i, m_i \in \Z_{>0}$. 
Since $\psi((\xi_1, \cdots, \xi_r))=0$, 
for any $i, j \in \{1, \cdots, r\}$, an equation  
\begin{equation}\label{e-psi-relation}
\begin{split}
\,\,\,\,&(f_if_j)^{k_{ij}}((f_if_j)^{\ell_{ij}}+y_{ij})(f_i^{m_i}f_j^{m_j}+f_i^{m_i}x_j)f_j^{n_j}\mu_i \\
=&(f_if_j)^{k_{ij}}((f_if_j)^{\ell_{ij}}+y_{ij})(f_i^{m_i}f_j^{m_j}+f_j^{m_j}x_i)f_i^{n_i}\mu_j
\end{split}
\end{equation}
holds in $M$ for some $y_{ij} \in I$, $\ell_{ij}, k_{ij} \in \Z_{>0}$.

Let $B_1$ be the finitely generated $\Z$-subalgebra of $A$ generated by $\{f_i\} \cup \{g_i\} \cup \{a_i\} \cup \{x_i\} \cup \{y_{ij}\}$. 
Let 
$$J:=\sum_{i=1}^r B_1x_i+\sum_{1 \leq i, j \leq r} B_1y_{ij} \subset B_1.$$
Let $B:=(1+J)^{-1}B_1$. 
We have a natural ring homomorphism $B \to A$. 
In particular, $M$ is an $B$-module. 
Let $N$ be the finitely generated $B$-submodule of $M$ generated by $\mu_1, \cdots, \mu_r$. 
Thanks to Step~\ref{s-noether}, the sequence 
$$0 \to N \xrightarrow{\varphi'} \prod^r_{i=1} N \otimes_B (1+JB_{f_i})^{-1}B_{f_i}
\xrightarrow{\psi'} \prod_{1 \leq i, j \leq r} N \otimes_B (1+JB_{f_if_j})^{-1}B_{f_if_j}$$
is exact. 
Take an element 
$$\eta:=\left(\frac{\frac{\mu_1}{f_1^{n_1}}}{1+\frac{x_1}{f_1^{m_1}}}, \cdots, \frac{\frac{\mu_r}{f_r^{n_r}}}{1+\frac{x_r}{f_i^{m_r}}}\right) \in 
\prod^r_{i=1} N \otimes_B (1+JB_{f_i})^{-1}B_{f_i}.$$
It follows from (\ref{e-psi-relation}) that $\psi'(\eta)=0$. 
Thus, we can find an element $\mu \in N$ such that $\varphi'(\mu)=\eta$. 
For any $i \in \{1, \cdots, r\}$, the equation 
$$f_i^{n_i+q_i}(f_i^{p_i}+z_i)(f_i^{m_i}+x_i)\mu=f_i^{m_i+q_i}(f_i^{p_i}+z_i)\mu_i$$
holds in $N$ for some $p_i, q_i \in \Z_{>0}$ and $z_i \in J$. 
Therefore, we get 
$$\varphi(\mu)=\left(\frac{\frac{\mu_1}{f_1^{n_1}}}{1+\frac{x_1}{f_1^{m_1}}}, \cdots, \frac{\frac{\mu_r}{f_r^{n_r}}}{1+\frac{x_r}{f_i^{m_r}}}\right),$$
as desired. 
This completes the proof of Step~\ref{s-general}. 
\end{proof}
Step~\ref{s-general} completes the proof of Theorem~\ref{t-FK-sheafy}. 
\end{proof}

\subsubsection{Global case}\label{ss-FK-global}

Let $X$ be a scheme and let $V$ be a closed subset of $X$ equipped with the induced topology from $X$. 
Set $i:V \hookrightarrow X$ to be the induced continuous map. 
For a quasi-coherent sheaf $F$ on $X$ and an affine open subset $U$ of $X$ with $A:=\Gamma(U, \MO_X)$, 
it follows from Theorem~\ref{t-FK-sheafy} that 
\begin{equation}\label{e-FK-sheafy}
\Gamma(V \cap U, i^{-1}(F))=\Gamma(U, F) \otimes_A (1+IA)^{-1}A, 
\tag{\ref{ss-FK-global}.1}
\end{equation}
where $I$ is an ideal of $A$ such that 
the closed subset $V(I)$ of $U=\Spec\,A$ corresponding to $I$ 
is equal to $V \cap U$.

\section{Zariskian f-adic rings}

In Subsection~\ref{ss-zar-def}, 
we introduce Zariskian f-adic rings and Zariskisation of f-adic rings. 
In Subsection \ref{ss-Zariskisation}, 
we show that the Zariskisation of any f-adic ring is Zariskian. 
In Subsection~\ref{ss-zar-haus} (resp. \ref{ss-zar-comp}), 
we discuss the relation between Zariskisations and 
Hausdorff quotients (resp. completions). 
In Subsection~\ref{ss-zar-aff}, 
we introduce Zariskian affinoid rings. 

For some corresponding results for the case of adic rings, 
we refer to \cite[Ch. 0, Section 7.3(b) and Ch. I, Section B.1]{FK}. 

\subsection{Definition}\label{ss-zar-def}

In this subsection, we introduce Zariskian f-adic rings and 
Zariskisation of f-adic rings.

\subsubsection{Definition as rings}\label{sss-zar-def1}

Let us first define Zariskisations of f-adic rings 
without topological structures. 
In (\ref{sss-zar-def2}), we shall introduce topologies on them. 

\begin{dfn}\label{d-zar}
Let $A$ be an f-adic ring. 
We set $S_A^{\Zar}:=1+A^{\circ\circ}$. 
Since $A^{\circ\circ}$ is an ideal of $A^{\circ}$, 
it holds that $S_A^{\Zar}$ is a multiplicative subset of $A$. 
We set $A^{\Zar}:=(S_A^{\Zar})^{-1}A$, 
which is called the {\em Zariskisation} of $A$. 
Also the natural ring homomorphism $A \to A^{\Zar}$ is 
called the {\em Zariskisation} of $A$. 
We say that $A$ is {\em Zariskian} if $S_A^{\Zar} \subset A^{\times}$. 
\end{dfn}

\begin{rem}\label{r-zar1}
Given an f-adic ring $A$, 
it is obvious that $A$ is Zariskian if and only if the natural ring homomorphism 
$A \to A^{\Zar}$ is bijective. 
\end{rem}

\begin{rem}\label{r-zar2}
If $\varphi:A \to B$ be a continuous ring homomorphism, 
then we have that $A^{\circ\circ} \subset B^{\circ\circ}$ 
and $\varphi(S_A^{\Zar}) \subset S_B^{\Zar}$. 
In particular, we get a commutative diagram of ring homomorphisms:
$$\begin{CD}
A @>\varphi >> B\\
@VVV @VVV\\
A^{\Zar} @>\varphi^{\Zar}>> B^{\Zar}, 
\end{CD}$$ 
where the vertical arrows are the natural ring homomorphisms. 
\end{rem}

\begin{lem}\label{l-zar-sub}
Let $A$ be an f-adic ring and let $A_0$ be an open subring of $A$. 
Then $A$ is Zariskian if and only if $A_0$ is Zariskian. 
\end{lem}

\begin{proof}
Assume that $A_0$ is Zariskian. 
Take $y \in A^{\circ\circ}$. 
We can find a positive integer $n$ such that $y^n \in A_0^{\circ\circ}$. 
Since $A_0$ is Zariskian, we have that 
$$1-y^n \in 1+A_0^{\circ\circ} \subset A_0^{\times} \subset A^{\times}.$$
This implies $1-y \in A^{\times}$. 
Thus $A$ is Zariskian. 

Conversely, assume that $A$ is Zariskian. 
Take $y \in A_0^{\circ\circ}$. 
We have that 
$$1+y \in 1+A_0^{\circ\circ} \subset 1+A^{\circ\circ} \subset A^{\times},$$
where the last inclusion holds because $A$ is Zariskian. 
Thus, there exists $w \in A$ such that $(1+y)w=1$. 
It suffices to show that $w \in A_0$. 
Let us prove $wy^m \in A_0$ for any $m \in \Z_{\geq 0}$ 
by descending induction on $m$. 
Since $y \in A_0^{\circ\circ}$, we get $wy^m \in A_0$ for $m \gg 0$. 
Thanks to the equation 
$$y^mw+y^{m+1}w=y^m \in A_0,$$
if $y^{m+1}w \in A_0$, then $y^mw \in A_0$. 
Therefore, it follows that $w \in A_0$. 
Thus $A_0$ is Zariskian. 
\end{proof}

\subsubsection{Definition as f-adic rings}\label{sss-zar-def2}

For any f-adic ring $A$, 
we introduce a topology on the ring $A^{\Zar}$ (Definition~\ref{d-zar-top}), 
so that also $A^{\Zar}$ is an f-adic ring. 
To this end, we need two auxiliary results: Lemma~\ref{l-two-S-inv} 
and Lemma~\ref{l-top-indep}. 
Furthermore, we show that $A^{\Zar}$ is 
an initial object of the category of Zariskian f-adic $A$-algebras 
(Theorem~\ref{t-zar-univ}).

\begin{lem}\label{l-two-S-inv}
Let $A$ be an f-adic ring. 
Let $P$ be an open pseudo-subring of $A^{\circ\circ}$. 
Then the natural ring homomorphism 
$\theta:(1+P)^{-1}A \to A^{\Zar}$ is bijective. 
\end{lem}

\begin{proof}
Take $z \in S_A^{\Zar}$. 
Since $S_A^{\Zar}=1+A^{\circ\circ}$, 
we can write $z=1-y$ for some $y \in A^{\circ\circ}$. 
It follows from $y \in A^{\circ\circ}$ that 
we get $y^k \in P$ for some positive integer $k$. 
Therefore, for any $a \in A$, we get 
$$\theta\left(\frac{a(1+y+\cdots+y^{k-1})}{1-y^k}\right)=\frac{a}{1-y}=\frac{a}{z}.$$
Thus $\theta$ is surjective. 

Take $a \in A$ and $w \in 1+P$ such that $\theta(a/w)=0$. 
This implies that $za=0$ for some $z \in S_A^{\Zar}$. 
For $y \in A^{\circ\circ}$ and $k \in \Z_{>0}$ such that 
$z=1-y$ and $y^k \in P$, 
we have that 
$(1-y^k)a=0$ and $1-y^k \in 1+P$. 
Therefore the equation $a/w=0$ holds, hence $\theta$ is injective. 
\end{proof}

\begin{rem}\label{r-FK-def}
If $A$ is an f-adic ring 
whose topology coincides with the $I$-adic topology 
for some ideal $I$ of $A$, 
then Lemma~\ref{l-two-S-inv} implies 
the following assertions. 
\begin{enumerate}
\item 
The natural ring homomorphism 
$(1+I)^{-1}A \to A^{\Zar}$ is bijective. 
\item 
$A$ is Zariskian if and only if $1+I \subset A^{\times}$. 
\end{enumerate}
\end{rem}

\begin{rem}\label{r-two-S-inv}
For an f-adic ring $A$ and an open subring $A'$ of $A$, 
we can consider $(A')^{\Zar}$ as a subring of $A^{\Zar}$ 
by Lemma~\ref{l-two-S-inv}. 
\end{rem}

\begin{lem}\label{l-top-indep}
Let $A$ be an f-adic ring. 
For any $i \in \{1, 2\}$, 
let $A_i$ be a ring of definition of $A$ and 
let $I_i$ be an ideal of definition of $A_i$. 
Then, for any positive integer $n_1$, 
there exists a positive integer $n_2$ such that 
the inclusion 
$$I_2^{n_2}A_2^{\Zar} \subset I_1^{n_1}A_1^{\Zar}$$
holds as subsets of $A^{\Zar}$ (cf. Remark~\ref{r-two-S-inv}). 
\end{lem}

\begin{proof}
Replacing $(A_2, I_2)$ by $(A_1 \cdot A_2, I_2A_1 \cdot A_2)$ \cite[Corollary 1.3(i)]{Hub93}, 
we may assume that $A_1 \subset A_2$. 
Fix a positive integer $n_1$. 
There is a positive integer $n_2$ such that $I_2^{n_2}A_2 \subset I_1^{n_1}A_1$. 
Then we get 
$$I_2^{n_2}A_2^{\Zar}=I_2^{n_2}(S^{\Zar}_{A_2})^{-1}A_2=I_2^{n_2}(S^{\Zar}_{A_1})^{-1}A_2 \subset I_1^{n_1}(S^{\Zar}_{A_1})^{-1}A_1=I_1^{n_1}A_1^{\Zar},$$
where the second equation holds, 
since the inclusion $A_1 \subset A_2$ enables us to apply Lemma \ref{l-two-S-inv}. \end{proof}

\begin{dfn}\label{d-zar-top}
Let $A$ be an f-adic ring and let $A^{\Zar}$ be its Zariskisation. 
For a ring of definition $A_0$ of $A$ and 
an ideal of definition $I_0$ of $A_0$, 
we equip $A^{\Zar}$ with the group topology 
defined by $\{I_0^kA_0^{\Zar}\}_{k\in \Z_{>0}}$. 
Thanks to Lemma~\ref{l-top-indep}, 
this topology does not depend on the choice of $A_0$ and $I_0$. 
It is easy to check that $A^{\Zar}$ is an f-adic ring (cf. Lemma \ref{l-top-criterion}). 
\end{dfn}

\begin{rem}\label{r-zar-adic}
For any f-adic ring $A$, the natural ring homomorphism $A \to A^{\Zar}$ is adic. 
\end{rem}

\begin{lem}\label{l-zar-cont}
Let $\varphi:A \to B$ be a continuous ring homomorphism of f-adic rings and 
let $\varphi^{\Zar}:A^{\Zar} \to B^{\Zar}$ be the induced ring homomorphism 
(cf. Remark~\ref{r-zar2}). 
Then $\varphi^{\Zar}$ is continuous. 
Moreover, if $\varphi$ is adic, so is $\varphi^{\Zar}$. 
\end{lem}

\begin{proof}
Take a ring of definition $A_0$ (resp. $B_0$) of $A$ (resp. $B$) 
and an ideal of definition $I_0$ (resp. $J_0$) of $A_0$ (resp. $B_0$). 
Fix a positive integer $m$. 
Since $\varphi$ is continuous, we can find a positive integer $n$ such that 
$$I_0^n \subset \varphi^{-1}(J_0^m).$$
It follows from Lemma \ref{l-two-S-inv} that 
$A_0^{\Zar}=(1+I_0^n)^{-1}A_0$ and $B_0^{\Zar}=(1+J_0^m)^{-1}B_0$. 
In particular, any element $\zeta \in I_0^nA_0^{\Zar}$ can be written 
by $\zeta=(1+x)^{-1}x'$ for some $x, x' \in I_0^n$. 
Then we get 
\[
\varphi^{\Zar}(\zeta)=\varphi^{\Zar}\left(\frac{x'}{1+x} \right)=\frac{\varphi(x')}{1+\varphi(x)}
\in J_0^m\cdot (1+J_0^m)^{-1}B_0= J_0^mB_0^{\Zar}.
\]
Hence, it holds that $I_0^nA^{\Zar}_0 \subset (\varphi^{\Zar})^{-1}(J_0^mB^{\Zar}_0).$ 
Therefore, $\varphi^{\Zar}$ is continuous. 

If $\varphi$ is adic, then \cite[Corollary 1.9(ii)]{Hub93} 
implies that $\varphi^{\Zar}$ is adic. 
\end{proof}

\begin{thm}\label{t-zar-univ}
Let $A$ be an f-adic ring and 
let $\theta:A \to A^{\Zar}$ be the Zariskisation. 
Then, for any continuous ring homomorphism $\varphi:A \to B$ 
to a Zariskian f-adic ring $B$, 
there exists a unique continuous ring homomorphism 
$\psi:A^{\Zar} \to B$ such that $\varphi=\psi \circ \theta$. 
\end{thm}

\begin{proof}
The assertion follows from Remark~\ref{r-zar1} and Lemma~\ref{l-zar-cont}. 
\end{proof}

\subsection{Zariskisations are Zariskian}\label{ss-Zariskisation}

In this subsection, 
we show that the Zariskisation of any f-adic ring is Zariskian 
(Theorem~\ref{t-zar-zar}).

\begin{lem}\label{l-frac-indep}
Let $A$ be an f-adic ring. 
Take elements $a_1, a_2 \in A$ and $s_1, s_2 \in S^{\Zar}_A$. 
If $a_1 \in A^{\circ\circ}$ and 
the equation $a_1/s_1=a_2/s_2$ holds in $A^{\Zar}$, 
then $a_2 \in A^{\circ\circ}$. 
\end{lem}

\begin{proof}
We have the equation $s_3s_2a_1=s_3s_1a_2$ in $A$ for some $s_3 \in S^{\Zar}_A$. 
As we can write $s_i=1+y_i$ for some $y_i \in A^{\circ\circ}$, 
we get 
$$(1+y_3)(1+y_2)a_1=(1+y_3)(1+y_1)a_2.$$
Thanks to $a_1 \in A^{\circ\circ}$ and the fact that $A^{\circ\circ}$ is an ideal of $A^{\circ}$, 
we can find $y \in A^{\circ\circ}$ such that 
$$(1+y)a_2 \in A^{\circ\circ}.$$

We show that $y^na_2 \in A^{\circ\circ}$ for any non-negative integer $n$ 
by descending induction on $n$. 
If $n \gg 0$, then it follows from $y \in A^{\circ\circ}$ that 
$y^na_2 \in A^{\circ\circ}$. 
Assume that $y^{k+1}a_2 \in A^{\circ\circ}$ for a non-negative integer $k$.  
Since 
$$y^ka_2+y^{k+1}a_2=y^k(1+y)a_2 \in A^{\circ\circ},$$
we have that $y^ka_2 \in A^{\circ\circ}$. 
Therefore, we get $a_2 \in A^{\circ\circ}$, as desired. 
\end{proof}

\begin{prop}\label{p-zariski-calcu}
Let $A$ be an f-adic ring. 
Then the equation 
$$(A^{\Zar})^{\circ\circ} = (S^{\Zar}_A)^{-1}A^{\circ\circ}$$
holds, where $(S^{\Zar}_A)^{-1}A^{\circ\circ}:=\{s^{-1}y \in A^{\Zar}\,|\,y \in A^{\circ\circ}, s \in S^{\Zar}_A\}$. 
\end{prop}

\begin{proof}
First, we show $(A^{\Zar})^{\circ\circ} \supset (S^{\Zar}_A)^{-1}A^{\circ\circ}$. 
Take $z \in (S^{\Zar}_A)^{-1}A^{\circ\circ}$ 
and we prove $z \in (A^{\Zar})^{\circ\circ}$. 
We can write 
$$z=\frac{y_1}{1+y_2} \in A^{\Zar}$$ 
for some $y_1, y_2 \in A^{\circ\circ}$. 
Let $A_0$ be a ring of definition of $A$ satisfying $y_1, y_2 \in A_0$. 
We see that $y_1, y_2 \in A_0^{\circ\circ}$. 
Take an ideal of definition $I_0$ of $A_0$. 
After replacing $z$ by a sufficiently large power $z^N$, 
we may assume that $y_1 \in I_0$. 
In particular, we have that 
$$z \in I_0A_0^{\Zar} \subset (A^{\Zar})^{\circ\circ},$$
where the inclusion follows from the definition of the topology on $A^{\Zar}$ (cf. Definition \ref{d-zar-top}). 
Thus the inclusion $(A^{\Zar})^{\circ\circ} \supset (S^{\Zar}_A)^{-1}A^{\circ\circ}$ 
holds. 

Second, we show $(A^{\Zar})^{\circ\circ} \subset (S^{\Zar}_A)^{-1}A^{\circ\circ}$. 
Take $z \in (A^{\Zar})^{\circ\circ}$. 
We can write  $z=a/s$ for some $a \in A$ and $s \in S^{\Zar}_A$. 
Let $A_0$ be a ring of definition of $A$ and 
let $I_0$ be an ideal of definition of $A_0$. 
Since $z^n \in I_0A_0^{\Zar}$ for some $n \in \Z_{>0}$, 
we can write 
$$z^n=\sum_{\ell=1}^m \frac{y_{1\ell}}{1+y_{2\ell}}$$ 
for some $y_{1\ell} \in I_0, y_{2\ell} \in A_0^{\circ\circ}$. 
In particular, we get $y_{1\ell}, y_{2\ell} \in A^{\circ\circ}$. 
Hence, we obtain $a^n/s^n=z^n \in (S^{\Zar}_A)^{-1}A^{\circ\circ}$. 
It follows from Lemma~\ref{l-frac-indep} that $a^n \in A^{\circ\circ}$. 
Since $A^{\circ\circ}=\sqrt{A^{\circ\circ}}$, we get $a \in A^{\circ\circ}$. 
Therefore, we have that $z=a/s \in (S^{\Zar}_A)^{-1}A^{\circ\circ}$, as desired. 
\end{proof}

\begin{thm}\label{t-zar-zar}
Let $A$ be an f-adic ring. 
Then the Zariskisation $A^{\Zar}$ of $A$ is Zariskian. 
\end{thm}

\begin{proof}
Take $z \in 1+(A^{\Zar})^{\circ\circ}$. 
By Proposition~\ref{p-zariski-calcu}, we can write 
$$z=1+\frac{y_1}{1+y_2}=\frac{1+y_1+y_2}{1+y_2}$$ 
for some $y_1, y_2 \in A^{\circ\circ}$. 
Since $1+y_1+y_2 \in S_A^{\Zar}$, 
we have that $z \in (A^{\Zar})^{\times}$. 
Hence, $A^{\Zar}$ is Zariskian. 
\end{proof}

\subsection{Relation to Hausdorff quotient}\label{ss-zar-haus}

In this subsection, we discuss the relation between 
Hausdorff quotients and Zariskisations. 
It is easy to find a Hausdorff f-adic ring that is not Zariskian 
(e.g. the integer ring $\Z$ with the $p$-adic topology). 
On the other hand, if a Zariskian f-adic ring $A$ 
is noetherian and its topology coincides with the $I$-adic topology for some ideal $I$, 
then $A$ is Hausdorff (Example~\ref{e-sep-zar}(1)). 
Unfortunately the same statement is false in general 
(Example~\ref{e-sep-zar}(3)). 
On the other hand, these two operations commute each other (Proposition~\ref{p-hd-zar}).

\begin{lem}\label{l-zar-sep}
Let $A$ be an f-adic ring and let $\alpha:A \to A^{\Zar}$ 
be the natural ring homomorphism. 
Then the inclusion $\Ker(\alpha) \subset \overline{\{0\}}$ holds.
\end{lem}

\begin{proof}
Take $x \in \Ker(\alpha)$. 
Then $(1-y)x=0$ for some $y \in A^{\circ\circ}$, which implies 
$x=y^kx$ for any positive integer $k$. 
Hence, we get $x \in \overline{\{0\}}$. 
\end{proof}

\begin{ex}\label{e-sep-zar}
Let $A$ be a ring and let $I$ be an ideal of $A$. 
We equip $A$ with the $I$-adic topology. 
\begin{enumerate}
\item 
If $A$ is a noetherian ring, then it follows from \cite[Theorem 10.17]{AM69} 
that 
$$\bigcap_{n=1}^{\infty} I^n=\{ x \in A\,|\,(1+y)x=0\text{ for some }y\in I\}.$$
This equation implies that the kernels of 
$A \to A^{\hd}$ and $A \to (1+I)^{-1}A$ coincide. 
In particular, if $A$ is noetherian and Zariskian, 
then $A$ is Hausdorff. 
\item 
In general, the inclusion 
$$\bigcap_{n=1}^{\infty} I^n \supset \{ x \in A\,|\,(1+y)x=0\text{ for some }y\in I\}$$
holds by Lemma~\ref{l-zar-sep}. 
\item 
By \cite[Remark (2) after Theorem 10.17]{AM69}, 
there exist a ring $B$ and an ideal $J$ of $B$ such that 
$$\bigcap_{n=1}^{\infty} J^n \supsetneq \{ x \in B\,|\,(1+y)x=0\text{ for some }y\in J\}.$$
We equip $B$ with the $J$-adic topology. 
Then the Zariskisation $B^{\Zar}=(1+J)^{-1}B$ of $B$ 
is not Hausdorff. 
We shall further study this example in Subsection~\ref{ss-non-surje}. 
\end{enumerate}
\end{ex}

\begin{lem}\label{l-zar-quot}
Let $A$ be an f-adic ring and let $J$ be an ideal of $A$. 
Let $A/J$ be the quotient f-adic ring of $A$ by $J$ (cf. Subsection \ref{ss-quot-fad}). 
If $A$ is Zasikian, then also $A/J$ is Zariskian. 
\end{lem}

\begin{proof}
Let $\pi:A \to A/J$ be the induced ring homomorphism. 
Fix a ring of definition $A_0$ of $A$ 
and an ideal of definition $I_0$ of $A_0$. 
We have that 
$$1+\pi(I_0A_0)=\pi(1+I_0A_0) \subset \pi(1+A^{\circ\circ}) 
\subset \pi(A^{\times}) \subset (A/J)^{\times}.$$  
Therefore, $A/J$ is Zariskian by Lemma~\ref{l-two-S-inv}. 
\end{proof}

\begin{lem}\label{l-hd-zar}
Let $A$ be a Hausdorff f-adic ring. 
Then the following hold. 
\begin{enumerate}
\item The natural ring homomorphism $A \to A^{\Zar}$ is injective. 
\item $A^{\Zar}$ is Hausdorff. 
\end{enumerate}
\end{lem}

\begin{proof}
The assertion (1) holds by Lemma~\ref{l-zar-sep}. 
Let us show (2). 
Take $\zeta \in \bigcap_{n \in \Z_{>0}} I_0^n A_0^{\Zar}$. 
Since $\zeta \in A_0^{\Zar}$, 
we can write $\zeta=a/s$ for some $a \in A_0$ and $s \in S_{A_0}^{\Zar}$. 
It suffices to show that $a \in \bigcap_{n \in \Z_{>0}} I_0^n$. 
Fix $n \in \Z_{>0}$. 
Since $\zeta \in I_0^nA_0^{\Zar}$, we get 
$a/s=x_n/s_n$ for some $x_n \in I_0^n$ and $s_n \in 1+I_0$. 
In particular, we get an equation 
$$ts_na=tsx_n \in I_0^n$$
for some $t \in 1+I_0$. 
We can write $ts_n=1-y$ for some $y \in I_0$, hence 
$$(1-y^{n})a=(1+y+\cdots+y^{n-1})tsx_n \in I_0^n.$$
As $y^n a \in I_0^n$, we get $a \in I_0^n$. 
Thus (2) holds. 
\end{proof}

\begin{prop}\label{p-hd-zar}
Let $A$ be an f-adic ring. 
Then both the natural continuous ring homomorphisms 
$$\theta_1:A^{\hd} \otimes_A A^{\Zar} \to (A^{\hd})^{\Zar} 
\quad\text{and}\quad
\theta_2:A^{\hd} \otimes_A A^{\Zar} \to (A^{\Zar})^{\hd}$$
are isomorphisms of topological rings, 
where $A^{\hd} \otimes_A A^{\Zar}$ denotes the topological tensor product (cf. Theorem~\ref{t-top-tensor}). 
\end{prop}

\begin{proof}
Fix a ring of definition $A_0$ of $A$ and 
an ideal of definition $I_0$ of $A_0$. 

\setcounter{step}{0}

\begin{step}\label{s-theta1}
The map 
$\theta_1:A^{\hd} \otimes_A A^{\Zar} \to (A^{\hd})^{\Zar}$ is an isomorphism 
of topological rings. 
\end{step}

\begin{proof}(of Step~\ref{s-theta1}) 
Let $N:=\Ker(A \to A^{\hd})$. 
Since the image of $A^{\circ\circ}$ by $A \to A^{\hd}$ is equal to $(A^{\hd})^{\circ\circ}$, we get 
$$A^{\hd} \otimes_A A^{\Zar} =(A/N) \otimes_A (S_A^{\Zar})^{-1}A
= (S^{\Zar}_{A^{\hd}})^{-1}(A/N)=(A^{\hd})^{\Zar},$$
hence the ring homomorphism $\theta_1$ is bijective. 
For any $k \in \Z_{\geq 0}$, the images of 
$I_0^k(A_0^{\hd} \otimes_{A_0} A_0^{\Zar})$ and $I_0^k(A_0^{\hd})^{\Zar}$ 
in $(A^{\hd})^{\Zar}$ coincide via $\theta_1$. 
Therefore, $\theta_1$ is an open map. 
This completes the proof of Step~\ref{s-theta1}. 
\end{proof}

\begin{step}\label{s-zh-z}
$(A^{\Zar})^{\hd}$ is Zariskian. 
\end{step}

\begin{proof}(of Step~\ref{s-zh-z}) 
The assertion follows from Lemma~\ref{l-zar-quot}. 
\end{proof}

\begin{step}\label{s-initial}
The map $\theta_2:A^{\hd} \otimes_A A^{\Zar} \to (A^{\Zar})^{\hd}$ is an isomorphism of topological rings. 
\end{step}

\begin{proof}(of Step~\ref{s-initial}) 
Let $\mathcal R$ be the category of 
Hausdorff Zariskian f-adic $A$-algebras 
whose arrows are continuous $A$-algebra homomorphisms. 
It follows from Lemma \ref{l-hd-zar}(2) that $(A^{\hd})^{\Zar}$ is Hausdorff. 
In particular, Step~\ref{s-theta1} implies that 
$A^{\hd} \otimes_A A^{\Zar}$ is a Hausdorff Zariskian f-adic $A$-algebra. 
Then we see that $A^{\hd} \otimes_A A^{\Zar}$ is 
an initial object of $\mathcal R$ 
by Remark~\ref{r-hd-univ}, Theorem~\ref{t-top-tensor} 
and Theorem~\ref{t-zar-univ}. 
It suffices to show that $(A^{\Zar})^{\hd}$ is an initial object of $\mathcal R$. 
It follows from Step~\ref{s-zh-z} that 
$(A^{\Zar})^{\hd}$ is an object of $\mathcal R$. 
By Remark~\ref{r-hd-univ} and Theorem~\ref{t-zar-univ}, 
$(A^{\Zar})^{\hd}$ is an initial object of $\mathcal R$, as desired. 
\end{proof}
Step~\ref{s-theta1} and Step~\ref{s-initial} complete 
the proof of Proposition~\ref{p-hd-zar}. 
\end{proof}

\subsection{Relation to completion}\label{ss-zar-comp}

The main result of this subsection 
is Theorem~\ref{t-comp-factor} 
which asserts that any f-adic ring and its Zariskisation 
has the same completion. 
We start with the following basic result.

\begin{lem}\label{l-comp-zar}
Let $A$ be a complete f-adic ring. 
Then $A$ is Zariskian. 
\end{lem}

\begin{proof}
Take $y \in A^{\circ\circ}$. 
Since $1+y+y^2+\cdots \in \widehat{A}$, we get $1-y \in A^{\times}$.  
\end{proof}

\begin{thm}\label{t-comp-factor}
Let $A$ be an f-adic ring. 
Let $\alpha:A \to A^{\Zar}$ be its Zariskisation and 
let $\gamma:A \to \widehat{A}$ be its completion. 
Then the following hold. 
\begin{enumerate}
\item 
There exists a unique continuous 
ring homomorphism $\beta:A^{\Zar} \to \widehat{A}$ such that $\gamma=\beta \circ \alpha$. 
\item 
For any non-negative integer $k$, 
ring of definition $A_0$ of $A$ and ideal of definition $I_0$ of $A_0$, 
the equation 
$$I_0^kA_0^{\Zar}=\beta^{-1}(I_0^k \widehat{A_0})$$
holds.   
\item 
The induced map $\widehat{\alpha}:\widehat{A} \to \widehat{A^{\Zar}}$ 
is an isomorphism of topological rings. 
\end{enumerate}
\end{thm}

\begin{proof}
The assertion (1) holds by Theorem~\ref{t-zar-univ} 
and Lemma~\ref{l-comp-zar}. 

We now show (2). 
Taking the Hausdorff quotients of $A \to A^{\Zar}$, 
we may assume that $A$ is Hausdorff and 
both $A$ and $A^{\Zar}$ are subrings of $\widehat{A}$ (Proposition~\ref{p-hd-zar}). 
It suffices to show the equation: 
$I_0^kA_0^{\Zar}=I_0^k \widehat{A_0} \cap A^{\Zar}.$ 
Since the inclusion $I_0^kA_0^{\Zar} \subset I_0^k \widehat{A_0} \cap A^{\Zar}$ is clear, 
it is enough to prove the opposite one: 
$I_0^kA_0^{\Zar} \supset I_0^k \widehat{A_0} \cap A^{\Zar}$. 

We have the following: 
$$I_0^k = I_0^k \widehat{A_0} \cap A \subset A \subset \widehat{A}.$$ 
Applying the functor $(-) \otimes_{A_0} A_0^{\Zar}$, we get  
{\small 
\begin{equation}\label{t-comp-factor1}
I_0^k \otimes_{A_0} A_0^{\Zar} = 
(I_0^k\widehat{A_0} \cap A) \otimes_{A_0} A_0^{\Zar} \subset A \otimes_{A_0} A_0^{\Zar} \subset \widehat{A} \otimes_{A_0} A_0^{\Zar}.
\end{equation}
}
Consider the induced ring isomorphism: 
$$\theta:\widehat{A} \otimes_{A_0} A_0^{\Zar} \xrightarrow{\simeq} \widehat{A}, \quad x \otimes y \mapsto xy.$$

\begin{claim}
The inclusion 
$$I_0^k\widehat{A_0} \cap A^{\Zar} \subset 
\theta((I_0^k\widehat{A_0} \cap A) \otimes_{A_0} A_0^{\Zar}).$$
holds. 
\end{claim}

\begin{proof}(of Claim) 
Take $\zeta \in I_0^k\widehat{A_0} \cap A^{\Zar}$. 
By Lemma~\ref{l-two-S-inv}, 
we can find $s \in 1+I_0$ such that $s\zeta \in I_0^k\widehat{A_0} \cap A$. 
In particular, we get 
$$s\zeta \otimes s^{-1} \in (I_0^k\widehat{A_0} \cap A) \otimes_{A_0} A_0^{\Zar},$$
which implies 
$$\zeta= \theta(s\zeta \otimes s^{-1}) \in \theta((I_0^k\widehat{A_0} \cap A) \otimes_{A_0} A_0^{\Zar}).$$
This completes the proof of Claim. 
\end{proof}

We obtain 
$$I_0^k\widehat{A_0} \cap A^{\Zar} \subset 
\theta((I_0^k\widehat{A_0} \cap A) \otimes_{A_0} A_0^{\Zar})
=\theta(I_0^k \otimes_{A_0} A_0^{\Zar})=I^k_0A_0^{\Zar},$$
where the inclusion holds by Claim and 
the first equation follows from (\ref{t-comp-factor1}). 
Thus (2) holds.

Thanks to (1) and (2), we can apply Lemma~\ref{l-common-comp}, 
hence the assertion (3) holds. 
This completes the proof of Theorem~\ref{t-comp-factor}. 
\end{proof}

\begin{cor}\label{c-zar-bdd-open}
Let $A$ be an f-adic ring. 
Let $\alpha:A \to A^{\Zar}$ be the Zariskisation of $A$. 
Then the following hold. 
\begin{enumerate}
\item 
The induced map $\alpha^*:\mathfrak O_{A^{\Zar}} \to \mathfrak O_{A}$ is bijective, and the inverse map $(\alpha^*)^{-1}=:\alpha_*$ satisfies $\alpha_*(A_0)=A_0^{\Zar}$ for any $A_0 \in \mathfrak O_{A}$. 
\item 
It holds that 
$\alpha_*(\mathfrak B_A)=\mathfrak B_{A^{\Zar}}$ 
and 
$\alpha^*(\mathfrak B_{A^{\Zar}})=\mathfrak B_{A}$. 
\item 
It holds that 
$\alpha_*(\mathfrak I_{A})=\mathfrak I_{A^{\Zar}}$ 
and 
$\alpha^*(\mathfrak I_{A^{\Zar}})=\mathfrak I_{A}$. 
\end{enumerate}
\end{cor}

\begin{proof}
Take the Zariskisation and the completion of $A$ (cf. Theorem~\ref{t-comp-factor}(1)): 
$$\gamma:A \xrightarrow{\alpha} A^{\Zar} \xrightarrow{\beta} \widehat{A}.$$

We show (1).   
By Lemma~\ref{l-complete-open}, 
we see that $\beta^*$ and $\gamma^*$ are bijective, hence so is $\alpha^*$. 
Take $A_0 \in \mathfrak O_{A}$. 
By Theorem~\ref{t-comp-factor}(3), 
we have that $\gamma_*(A_0)=\beta_*(A_0^{\Zar})$, 
which implies that 
$$\alpha_*(A_0)=\beta^*\beta_*\alpha_*(A_0)=\beta^*\gamma_*(A_0)=\beta^*\beta_*(A_0^{\Zar})=A_0^{\Zar}.$$
Thus (1) holds.

The assertion (2) (resp. (3)) follows from Lemma \ref{l-complete-open2} 
(resp. Lemma \ref{l-complete-open3}). 
\end{proof}

\subsection{Zariskian affinoid rings}\label{ss-zar-aff}

In this subsection, 
we introduce Zariskian affinoid rings and Zariskisation of affinoid rings. 

\begin{dfn}\label{d-zar-affinoid}
Let $A=(A^{\rhd}, A^+)$ be an affinoid ring. 
We define the {\em Zariskisation} $A^{\Zar}$ of $A$ 
by $A^{\Zar}:=((A^{\rhd})^{\Zar}, (A^+)^{\Zar})$. 
\end{dfn}

\begin{rem}\label{r-aff-zar-c}
We use the same notation as in Definition~\ref{d-zar-affinoid}. 
Then the following hold. 
\begin{enumerate}
\item 
By Definition~\ref{d-zar-top}, 
$(A^{\rhd})^{\Zar}$ is an f-adic ring. 
It follows from Corollary~\ref{c-zar-bdd-open}(3) that 
$(A^+)^{\Zar} \in \mathfrak I_{A^{\Zar}}$. 
Thus 
$$A^{\Zar}=((A^{\rhd})^{\Zar}, (A^+)^{\Zar})$$ 
is an affinoid ring. 
\item 
By Theorem~\ref{t-comp-factor}, there exist natural continuous $A$-algebra homomorphisms: 
$$\gamma:A \xrightarrow{\alpha} A^{\Zar} \xrightarrow{\beta} \widehat{A},$$
where $\beta$ and $\gamma$ are the completions. 
Moreover, all of $\alpha$, $\beta$ and $\gamma$ are adic. 
\end{enumerate}
\end{rem}

\begin{dfn}
An affinoid ring $A$ is {\em Zariskian} if 
the natural homomorphism $A \to A^{\Zar}$ is an isomorphism 
of affinoid rings  
(cf. Remark~\ref{r-aff-zar-c}(2)). 
\end{dfn}

\begin{rem}\label{r-zar-zar}
Let $A$ be an affinoid ring. 
Then its Zariskisation $A^{\Zar}$ is 
Zariskian by Theorem~\ref{t-zar-zar}. 
\end{rem}

\begin{prop}\label{p-zar-univ}
Let $A$ be an affinoid ring and 
let $\alpha:A \to A^{\Zar}$ be the Zariskisation of $A$. 
Then, for a continuous ring homomorphism $\varphi:A \to B$ 
to a Zariskian affinoid ring $B$, 
there  exists a unique continuous ring homomorphism 
$\psi:A^{\Zar} \to B$ of affinoid rings such that 
$\varphi=\psi \circ \alpha$. 
\end{prop}

\begin{proof}
The assertion follows from Theorem~\ref{t-zar-univ}. 
\end{proof}

\begin{lem}\label{l-zar-quot2}
Let $A=(A^{\rhd}, A^+)$ be an f-adic ring and let $J^{\rhd}$ 
be an ideal of $A^{\rhd}$. 
Let $A/J^{\rhd}$ be the quotient affinoid ring of $A$ by $J^{\rhd}$ (cf. Subsection \ref{ss-quot-aff}). 
If $A$ is Zariskian, then also $A/J^{\rhd}$ is Zariskian. 
\end{lem}

\begin{proof}
The assertion immediately follows from Lemma~\ref{l-zar-quot}. 
\end{proof}

\section{Zariskian adic spaces}

In this section, we introduce Zariskian adic spaces. 
In Subsection~\ref{ss-affinoid}, 
we establish foundations for the affinoid case. 
In particular we define a presheaf $\MO_A^{\Zar}$ on $\Spa\,A$ 
for any affinoid ring $A$. 
In Subsection~\ref{ss-structure}, 
we prove that $\MO_A^{\Zar}$ is actually a sheaf. 
In Subsection~\ref{ss-zar-ad-sp}, 
we define Zariskian adic spaces. 
In Subsection~\ref{ss-cex-Tate}, 
we exhibit examples that violate the Tate acyclicity. 

\subsection{Affinoid case}\label{ss-affinoid}

In this subsection, we introduce a presheaf $\MO_A^{\Zar}$ on $\Spa\,A$ 
for any affinoid ring $A$. 
This is a Zariskian analogue of Huber's affinoid adic spaces 
introduced in \cite[Section 1]{Hub94}. 
Indeed, most of our arguments are quite similar 
to the ones in \cite[Section 1]{Hub94}. 

In (\ref{sss-rat-localise}), we consider 
what the definition of $\MO_A^{\Zar}(U)$ should be for rational subsets $U$. 
In (\ref{sss-zar-str}), we define $\MO_A^{\Zar}$  
based on results obtained in (\ref{sss-rat-localise}).

\subsubsection{Rational localisation}\label{sss-rat-localise}

Let $A=(A^{\rhd}, A^+)$ be an affinoid ring and 
let $U$ be a rational subset of $\Spa\,A$. 
Then there exist elements $f_1, \cdots, f_r, g \in A^{\rhd}$ such that 
$(f_1, \cdots, f_r, g)$ is an open ideal of $A^{\rhd}$ and 
$$U=R\left(\frac{f_1, \cdots, f_r}{g}\right).$$
Take a ring of definition $A_0$ of $A^{\rhd}$ 
and an ideal of definition $I_0$ of $A_0$. 
We set $B^{\rhd}:=A^{\rhd}[1/g]=(A^{\rhd})_g$ and  
$B_0:=A_0\left[f_1/g, \cdots, f_r/g\right] \subset B^{\rhd}.$ 
We equip $B^{\rhd}$ with the group topology induced by $\{I^k_0 B_0\}_{k \in \Z_{>0}}$. 
We can check that this topology does not depend on the choice 
of $A_0$ and $I_0$. 
Set $B^+$ to be the integral closure of 
$B^+[f_1/g, \cdots, f_r/g]$ in $B^{\rhd}$. 
We set 
$$A\left(\frac{f_1, \cdots, f_r}{g}\right):=(B^{\rhd}, B^+).$$
It follows that $A\left(\frac{f_1, \cdots, f_r}{g}\right)$ is 
an affinoid ring and the induced ring homomorphism 
$A \to A\left(\frac{f_1, \cdots, f_r}{g}\right)$ 
is adic. 
Let $A\left(\frac{f_1, \cdots, f_r}{g}\right)^{\Zar}$ be 
the Zariskisation of $A\left(\frac{f_1, \cdots, f_r}{g}\right)$ and 
let 
$$\varphi=\varphi_U:A \xrightarrow{\alpha} A\left(\frac{f_1, \cdots, f_r}{g}\right) 
\xrightarrow{\theta} A\left(\frac{f_1, \cdots, f_r}{g}\right)^{\Zar}=:Z_A(U)$$
be the induced adic ring homomorphisms. 
Although it is not clear that $Z_A(U)$ does not depend on the choices 
of $f_1, \cdots, f_r, g$, 
we shall later see in Lemma~\ref{l-rat-universal1} that this is actually true. 
To this end, we first establish an auxiliary result: 
Lemma~\ref{l-zar-nonempty}.

\begin{lem}\label{l-zar-nonempty}
Let $A=(A^{\rhd}, A^+)$ be a Zariskian affinoid ring. 
Then the following hold. 
\begin{enumerate}
\item $(A^{\rhd})^{\times}$ is an open subset of $A^{\rhd}$. 
\item If $\m$ is a maximal ideal of $A^{\rhd}$, 
then there exists a point $v \in \Spa\,A$ such that $\m=\Ker(v)$. 
\end{enumerate}
\end{lem}

\begin{proof}
First we show (1). 
Since $1+(A^{\rhd})^{\circ\circ}$ is an open subset of $A^{\rhd}$, 
so is 
$$(A^{\rhd})^{\times}=\bigcup_{x \in (A^{\rhd})^{\times}} x \cdot (1+(A^{\rhd})^{\circ\circ}).$$ 
Thus the assertion (1) holds. 

Second we show (2). 
Since $(A^{\rhd})^{\times}$ is an open subset of $A^{\rhd}$ by (1),  
the inclusion 
$\m \subset A^{\rhd} \setminus (A^{\rhd})^{\times}$ 
implies 
$\overline{\m} \subset A^{\rhd} \setminus (A^{\rhd})^{\times}$, 
where $\overline{\m}$ denotes the closure of $\m$ in $A^{\rhd}$. 
Since $\overline{\m}$ is again an ideal of $A^{\rhd}$ 
such that $\m \subset \overline{\m} \subsetneq A^{\rhd}$, 
we have that $\m=\overline{\m}$. 
We set $B^{\rhd}:=A^{\rhd}/\m$ to be the topological ring 
equipped with the quotient topology induced from $A^{\rhd}$. 
Since $\m=\overline{\m}$, it holds that $B^{\rhd}$ is Hausdorff. 
Let $B^+$ be the integral closure of the image of $A^+$. 
Then $B:= (B^{\rhd}, B^+)$ is an affinoid ring 
(cf. Subsection \ref{ss-quot-aff}). 
It follows from \cite[Proposition 3.6(i)]{Hub93} that 
$\Spa\,B \neq \emptyset$. 
Since we have a natural continuous ring homomorphism $A \to B$, 
any element of the image of $\Spa\,B \to \Spa\,A$ 
is an element $v \in \Spa\,A$ as in the statement. 
Thus (2) holds. 
\end{proof}

\begin{lem}\label{l-rat-universal1}
Let $A=(A^{\rhd}, A^+)$ be an affinoid ring. 
Let $U$ be a rational subset of $\Spa\,A$ 
and let $f_1, \cdots, f_r, g \in A^{\rhd}$ be elements 
such that $(f_1, \cdots, f_r, g)$ is an open ideal of $A^{\rhd}$ and 
$$U=R\left(\frac{f_1, \cdots, f_r}{g}\right).$$
If $\psi:A \to C$ is a continuous ring homomorphism 
to a Zariskian affinoid ring 
such that the image of the induced map $\Spa\,(\psi):\Spa\,C \to \Spa\,A$ 
is contained in $U$,  
then there exists a unique continuous ring homomorphism 
$$\psi':A\left(\frac{f_1, \cdots, f_r}{g}\right)^{\Zar} \to C$$ 
such that $\psi=\psi' \circ \varphi$.  
\end{lem}

\begin{proof}
Let $\psi:A \to C$ be a continuous ring homomorphism 
to a Zariskian affinoid ring $C=(C^{\rhd}, C^+)$ 
satisfying ${\rm Im}\,(\Spa\,(\psi)) \subset U$. 
We now show the following two assertions. 
\begin{enumerate}
\item[(i)] $\psi(g) \in (C^{\rhd})^{\times}$. 
\item[(ii)] $\psi(f_1)/\psi(g), \cdots, \psi(f_r)/\psi(g) \in C^+$. 
\end{enumerate}

First we prove (i). 
By $\Image(\Spa\,(\psi)) \subset U=R\left(\frac{f_1, \cdots, f_r}{g}\right)$, 
we have that $v(\psi(g)) \neq 0$ for any $v \in \Spa\,C$. 
It follows from Lemma~\ref{l-zar-nonempty}(2) that 
$\psi(g)$ is not contained in any maximal ideal of $C^{\rhd}$. 
Thus (i) holds. 

Second we show (ii). 
Again by $\Image(\Spa\,(\psi)) \subset U=R\left(\frac{f_1, \cdots, f_r}{g}\right)$, 
we have that $v(\psi(f_i)/\psi(g)) \leq 1$ 
for any $v \in \Spa\,C$ and any $i \in \{1, \cdots, r\}$. 
By \cite[Lemma 3.3(i)]{Hub93}, we get $\psi(f_i)/\psi(g) \in C^+$ 
for any $i \in \{1, \cdots, r\}$. 
Thus (ii) holds. 

\medskip

By (i) and (ii), 
there exists a unique ring homomorphism 
$\psi_1:A\left(\frac{f_1, \cdots, f_r}{g}\right) \to C$ 
inducing the following factorisation:  
$$\psi:A \xrightarrow{\alpha} A\left(\frac{f_1, \cdots, f_r}{g}\right) \xrightarrow{\psi_1} C.$$
It follows from \cite[(1.2)(ii)]{Hub94} that $\psi_1$ is continuous. 
Taking the Zariskisation of $\psi_1$, we get a unique factorisation of $\psi_1$: 
$$\psi:A \to A\left(\frac{f_1, \cdots, f_r}{g}\right) \xrightarrow{\theta} A\left(\frac{f_1, \cdots, f_r}{g}\right)^{\Zar} \xrightarrow{\psi'} C.$$
We are done. 
\end{proof}

\begin{lem}\label{l-rat-universal2}
Let $A$ be an affinoid ring and let $U$ be a rational subset of $\Spa\,A$. 
Then the following assertions hold. 
\begin{enumerate}
\item 
If $V$ is a rational subset of $\Spa\,A$ such that $V \subset U$, then 
there is a unique continuous ring homomorphism 
$Z_A(U) \to Z_A(V)$ that commutes with 
$\varphi_V:A \to Z_A(V)$ and $\varphi_U:A \to Z_A(U)$. 
\item 
The induced map $\Spa\,(\varphi_U):\Spa\,Z_A(U) \to \Spa\,A$ 
is an open injective map whose image is equal to $U$. 
If $W$ is a rational subset of $\Spa\,Z_A(U)$, then 
so is $\Spa\,(\varphi_U)(W)$. 
If $V$ is a rational subset of $\Spa\,A$ contained in $U$, then 
also $(\Spa\,(\varphi_U))^{-1}(V)$ is a rational subset. 
\item 
Set $B:=Z_A(U)$ and $g:=\Spa(\varphi_U):\Spa\,B \to \Spa\,A$. 
Let $V$ be a rational subset of $\Spa\,A$ contained in $U$. 
Then there exists a unique continuous ring homomorphism 
$r:Z_A(V) \to Z_B(g^{-1}(V))$ such that the following diagram is commutative: 
$$\begin{CD}
Z_B(g^{-1}(V)) @<<< Z_A(V)\\
@AAA @AAA\\
B @<<< A.
\end{CD}$$
Furthermore, $r$ is an isomorphism of topological rings. 
\end{enumerate}
\end{lem}

\begin{proof}
The assertion (1) holds by Lemma~\ref{l-rat-universal1}. 
The assertion (2) follows from \cite[Lemma 1.5(ii)]{Hub94} 
and Theorem \ref{t-comp-factor}. 
As for (3), 
we can apply the same argument as in \cite[Lemma 1.5(iii)]{Hub94} 
after replacing $F_A(-)$ and \cite[Lemma 1.3]{Hub94} by $Z_A(-)$ and 
Lemma~\ref{l-rat-universal1}, respectively. 
\end{proof}

\subsubsection{Zariskian structure presheaves}\label{sss-zar-str}

Let $A=(A^{\rhd}, A^+)$ be an affinoid ring. 
For any open subset $V$ of $\Spa\,A$, we set 
$$\Gamma(V, \MO_A^{\Zar}):=\MO_A^{\Zar}(V):=
\varprojlim_{\substack{U \subset V,\\ 
U:\text{ rational}\\ 
\text{subset}}} Z_A(U)^{\rhd},$$
where the inverse limit is taken in the category of $A^{\rhd}$-algebras. 
We equip $\MO_A^{\Zar}(V)$ with the inverse limit topology. 
If $V_1$ and $V_2$ are open subsets of $\Spa\,A$ satisfying 
$V_1 \supset V_2$, 
then the induced ring homomorphism  $\MO_A^{\Zar}(V_1) \to \MO_A^{\Zar}(V_2)$ is continuous. 
Thus $\MO^{\Zar}_A$ is a presheaf of topological rings. 

Fix $x \in \Spa\,A$. 
Set 
$$\MO^{\Zar}_{A, x}:=\varinjlim_{x \in U} \MO^{\Zar}_A(U),$$ 
where the direct limit is taken in the category of rings. 
We obtain a natural isomorphism of rings: 
$$\varinjlim_{\substack{x \in U,\\ 
U:\text{ rational}\\ 
\text{subset}}} \MO^{\Zar}_A(U) \xrightarrow{\simeq} 
\varinjlim_{x \in U} \MO^{\Zar}_A(U).$$
For any rational subset $U$ with $x \in U$, 
the valuation $x:A^{\rhd} \to \Gamma_x \cup \{0\}$ 
is uniquely extended to a valuation $v_U:\MO_A^{\Zar}(U) \to \Gamma_x \cup \{0\}$. 
Thus the set of valuations $\{v_U\}_{x \in U}$ define a valuation 
$$v_x:\MO_{A, x} \to \Gamma_x \cup \{0\}.$$
For any open subset $U$ of $\Spa\,A$, we set 
$$\MO_A^{\Zar, +}(U):=\{f \in \MO^{\Zar}_A(U)\,|\,v_x(f) \leq 1\text{ for any }x \in U\}.$$
Then $\MO_A^{\Zar, +}$ is a presheaf of rings on $\Spa\,A$. 
For any $x \in \Spa\,A$, 
we set $\MO_{A, x}^{\Zar, +}(U)$ to be the stalk of $\MO_A^{\Zar, +}$ at $x$.

\begin{prop}\label{p-zar-LRS}
The following hold. 
\begin{enumerate}
\item 
For any $x \in \Spa\,A$, 
the stalk $\MO_{A, x}^{\Zar}$ is a local ring 
whose maximal ideal is equal to $\Ker(v_x)$. 
\item 
For any $x \in \Spa\,A$, 
the stalk $\MO_{A, x}^{\Zar, +}$ is a local ring. 
It holds that 
$\MO_{A, x}^{\Zar, +}=\{f \in \MO_{A, x}^{\Zar}\,|\, v_x(f) \leq 1\}$ and 
its maximal ideal is equal to $\{f \in \MO_{A, x}^{\Zar}\,|\, v_x(f) < 1\}$. 
\item 
For any open subset $U$ of $\Spa\,A$ and $f, g \in \MO_A^{\Zar}(U)$, 
the set $\{x \in U\,|\,v_x(f) \leq v_x(g) \neq 0\}$ is 
an open subset of $\Spa\,A$. 
\item 
 For any rational subset $U$ of $\Spa\,A$, 
 it holds that $\MO_A^{\Zar}(U)=Z_A(U)^{\rhd}$ and 
$\MO_A^{\Zar, +}(U)=Z_A(U)^+$. 
\end{enumerate}
\end{prop}

\begin{proof}
We can apply the same argument as in \cite[Proposition 1.6]{Hub94} 
after replacing \cite[Lemma 1.5]{Hub94} by Lemma~\ref{l-rat-universal2}. 
\end{proof}

Let $A=(A^{\rhd}, A^+)$ be an affinoid ring. 
Let $M$ be an $A^{\rhd}$-module. 
For any open subset $V$ of $\Spa\,A$, we set 
$$\Gamma(V, M \otimes \MO_A^{\Zar}):=
\varprojlim_{\substack{U \subset V,\\ 
U:\text{ rational}\\ 
\text{subset}}}M \otimes_A Z_A(U)^{\rhd},$$
where the inverse limit is taken in the category of $A^{\rhd}$-modules. 
We have that $M \otimes \MO_A^{\Zar}$ is a presheaf of $\MO_A$-modules.

\subsection{Structure sheaves}\label{ss-structure}

The purpose of this subsection is to show that the structure presheaf $\MO_A^{\Zar}$ is actually a sheaf (Theorem~\ref{t-sheafy}). 
To this end, we start with an auxiliary result (Lemma \ref{l-rat-covering}) 
that assures the existence of refined rational coverings.

\begin{lem}\label{l-rat-covering}
Let $A=(A^{\rhd}, A^+)$ be a Zariskian affinoid ring and 
let $\mathcal U$ be an open cover of $\Spa\,A$. 
Then there exist elements $f_0, \cdots, f_n \in A^{\rhd}$ 
such that $\sum_{i=0}^n A^{\rhd} f_i=A^{\rhd}$ and 
the induced open cover $\left\{R\left(\frac{f_0, \cdots, f_n}{f_i}\right)\right\}_{0 \leq i \leq n}$ of $\Spa\,A$ is a refinement of $\mathcal U$. 
\end{lem}

\begin{proof}
Let $\gamma:A^{\rhd} \to \widehat{A^{\rhd}}$ be the completion. 
Thanks to \cite[Lemma 2.6]{Hub94}, we can find a refinement 
$\mathcal U$ by a rational covering 
$\left\{R\left(\frac{f_0, \cdots, f_n}{f_i}\right)\right\}_{0 \leq i \leq n}$ 
for some elements $f_0, \cdots, f_n, a_0, \cdots, a_n \in \widehat{A^{\rhd}}$ 
such that $\sum_{i=0}^n a_if_i=1$. 
Fix a ring of definition $A_0$ of $A^{\rhd}$ and an ideal of definition $I_0$ of $A_0$. 
By \cite[Lemma 3.10]{Hub93}, 
we can find elements $f'_0, \cdots, f'_n \in A^{\rhd}$ 
whose images by $\gamma$ are sufficiently close to $f_0, \cdots, f_n$, 
so that 
$$-1+\sum_{i=0}^na_i\gamma(f'_i) \in I_0\widehat{A_0}$$ 
and 
$$R\left(\frac{f_0, \cdots, f_n}{f_i}\right)=
R\left(\frac{f'_0, \cdots, f'_n}{f'_i}\right)$$
for any $i \in \{0, 1, \cdots, n\}$. 
Take elements $a'_i \in A^{\rhd}$ whose images by $\gamma$ 
are sufficiently close to $a_i$, so that 
$$-1+\sum_{i=0}^n\gamma(a'_if'_i) \in I_0\widehat{A_0}.$$ 
Therefore, we get 
$$-1+\sum_{i=0}^na'_if'_i \in \gamma^{-1}(I_0\widehat{A_0})=I_0A_0.$$
Since $A^{\rhd}$ is Zariskian, we get $\sum_{i=0}^nA^{\rhd} f'_i=A^{\rhd}$ as desired. 
\end{proof}

\begin{thm}\label{t-sheafy}
Let $A=(A^{\rhd}, A^+)$ be an affinoid ring and let $M$ be an $A^{\rhd}$-module. 
Then the presheaf $M \otimes \MO_A^{\Zar}$ on $\Spa\,A$ is a sheaf. 
\end{thm}

\begin{proof}
Replacing $A$ by its Zariskisation, 
we may assume that $A$ is Zariskian. 
Set $X:=\Spa\,A$. 
Take elements $f_0, \cdots, f_n \in A^{\rhd}$ such that 
$\sum_{k=0}^n A^{\rhd} f_k=A^{\rhd}$. 
Set $U_k:=R\left(\frac{f_0, \cdots, f_n}{f_k}\right)$ and 
$U_{k_1k_2} :=U_{k_1} \cap U_{k_2}$ for any $0 \leq k, k_1, k_2 \leq n$. 
Thanks to Lemma~\ref{l-rat-covering}, 
it suffices to show that the sequence 
{\small
\begin{equation}\label{e-sheafy}
0 \to (M \otimes \MO_A^{\Zar})(X) \xrightarrow{\varphi} \prod_{0\leq k \leq n} (M \otimes \MO_A^{\Zar})(U_k) 
\xrightarrow{\psi} \prod_{0\leq k_1<k_2 \leq n} (M \otimes \MO_A^{\Zar})(U_{k_1k_2})
\end{equation}
}
is exact. 
Fix a ring of definition $A_0$ of $A^{\rhd}$ and 
an ideal of definition $I_0$ of $A_0$.
We have that 
$$(M \otimes \MO_A^{\Zar})(X)=M,$$
$$(M \otimes \MO_A^{\Zar})(U_k)
=M \otimes_{A^{\rhd}} \left(1+I_0A_0\middle[\frac{f_0}{f_k}, \cdots, \frac{f_n}{f_k}\middle]\right)^{-1}A^{\rhd}\left[\frac{1}{f_k}\right],$$
$$(M \otimes \MO_A^{\Zar})(U_{k_1k_2})=M \otimes_{A^{\rhd}}
\left(1+I_0A_0\middle[\middle\{\frac{f_{\ell_1}f_{\ell_2}}{f_{k_1}f_{k_2}}\middle\}_{0 \leq \ell_1, \ell_2 \leq n}\middle]\right)^{-1}A^{\rhd}\left[\frac{1}{f_{k_1}f_{k_2}}\right].$$
Let $J_0$ be the $A_0$-submodule of $A^{\rhd}$ defined by 
$$J_0:=\sum_{k=0}^n A_0f_k \subset A^{\rhd}.$$
For $J:=J_0A^{\rhd}=A^{\rhd}$, we get a commutative diagram of schemes: 
$$\begin{CD}
\Proj(\bigoplus_{d=0}^{\infty} J^d) @>\pi >> \Spec\,A^{\rhd}\\
@VV\beta V @VV\alpha V\\
\Proj(\bigoplus_{d=0}^{\infty} J_0^d) @>\pi_0>> \Spec\,A_0,
\end{CD}$$
where $J^0:=A^{\rhd}$ and $J_0^0:=A_0$. 
By $J=A^{\rhd}$, we have that $\pi$ is an isomorphism. 
We get a morphism 
$$\sigma:=\beta \circ \pi^{-1}:\Spec\,A^{\rhd} \to 
\Proj\left(\bigoplus_{d=0}^{\infty} J_0^d\right).$$
It holds that $\beta$ is an affine morphism satisfying $\beta^{-1}(D_+(f_k))=D_+(f_k)$ for any $k \in \{0, \cdots, n\}$. 
Since the equation $\pi^{-1}(D(f_k))=D_+(f_k)$ holds, 
we obtain $\sigma^{-1}(D_+(f_k))=D(f_k)$. 
Thus, for any non-empty subset $K$ of $\{0, \cdots, n\}$, 
if we set 
$$f_K:=\prod_{k \in K} f_i,$$ 
then it holds that  
\begin{equation}\label{e-sheafy2}
\Gamma(D_+\left(f_K\right), \sigma_*\widetilde M)=
\Gamma(D\left(f_K\right), \widetilde M)=M \otimes_{A^{\rhd}} A^{\rhd}\left[\frac{1}{f_K}\right].
\end{equation}

\begin{claim}
For any non-empty subset $K$ of $\{0, \cdots, n\}$, 
the equation
\begin{equation}\label{e-sheafy3}
\Gamma\left(D_+(f_K), 
\MO_{\Proj(\bigoplus_{d=0}^{\infty} J_0^d)}\right)
=A_0\left[\middle\{\frac{f_{i_1}\cdots f_{i_{|K|}}}{f_K}\middle\}_{0 \leq i_j \leq n}\right] 
\end{equation}
holds, where the right hand side is defined as 
the $A_0$-subalgebra of $A^{\rhd}[1/f_K]$ generated by 
the set 
$\{f_{i_1}\cdots f_{i_{|K|}}f_K^{-1}\,|\, 0 \leq i_j \leq n\}$. 
\end{claim}

\begin{proof}(of Claim) 
Consider the following $A_0$-algebra homomorphisms of 
graded $A_0$-algebras preserving degrees 
\begin{eqnarray*}
A_0[t_0, \cdots, t_n] &\xrightarrow{\rho}& \bigoplus_{d=0}^{\infty} J_0^d 
\hookrightarrow A^{\rhd}[s]\\
t_k &\mapsto & f_ks
\end{eqnarray*}
where $A_0[t_0, \cdots, t_n]$ is the polynomial ring over $A_0$ and 
$f_ks$ is the element of degree one 
which is the same element as $f_k$. We set 
$$t_K:=\prod_{k \in K} t_k.$$
Taking localisations, we obtain $A_0$-algebra homomorphisms of 
graded $A_0$-algebras preserving degrees: 
\begin{eqnarray*}
A_0[t_0, \cdots, t_n, t_K^{-1}] &\xrightarrow{\rho_K}& 
\left(\bigoplus_{d=0}^{\infty} J_0^d\right)\left[\frac{1}{f_Ks^{|K|}}\right] 
\overset{q}\hookrightarrow \left(A^{\rhd}\left[\frac{1}{f_K}\right]\right)[s, s^{-1}].
\end{eqnarray*}
Since $\rho$ is surjective, so is $\rho_K$. 
It follows from \cite[Ch. II, Proposition 2.5(b)]{Har77} that 
the left hand side 
$\Gamma\left(D_+(f_K), \MO_{\Proj(\bigoplus_{d=0}^{\infty} J_0^d)}\right)$ 
of (\ref{e-sheafy3}) is nothing but the degree zero part of 
the graded ring $(\bigoplus_{d=0}^{\infty} J_0^d)[1/f_Ks^{|K|}]$, 
which is equal to 
$\rho_K\left(A_0[\{\frac{t_{i_1} \cdots t_{i_{|K|}}}{t_K}\}_{0 \leq i_j \leq n}]\right)$ 
because $\rho_K$ is surjective and preserving degrees. 
Embedding this ring via $q$, we get 
\begin{eqnarray*}
\Gamma\left(D_+(f_K), \MO_{\Proj(\bigoplus_{d=0}^{\infty} J_0^d)}\right)
&=&\rho_K\left(A_0\middle[\middle\{\frac{t_{i_1} \cdots t_{i_{|K|}}}{t_K}\middle\}_{0 \leq i_j \leq n}\middle]\right)\\
 &\xrightarrow{q,\, \simeq}& A_0\left[\middle\{\frac{f_{i_1} \cdots f_{i_{|K|}}}{f_K}\middle\}_{0 \leq i_j \leq n}\right]
\end{eqnarray*}
where the right hand side is the $A_0$-subalgebra of $A^{\rhd}[1/f_K]$ 
generated by $\{f_{i_1} \cdots f_{i_{|K|}}f_K^{-1}\,|\,0 \leq i_j \leq n\}$. 
This completes the proof of Claim. 
\end{proof}

Let 
$\iota: W=\pi^{-1}_0(V(I_0)) \hookrightarrow \Proj(\bigoplus_{d=0}^{\infty} J_0^d)$ be the inclusion map and we equip $W$ the induced topology. 
We set $\mathcal D(g):=D(g) \cap W$ for any homogeneous element $g$ 
of positive degree. 
For $\mathcal M:=\iota^{-1}\sigma_*\widetilde{M},$ 
we have that 
\begin{eqnarray*}
&&\Gamma(\mathcal D_+(f_k), \mathcal M)\\
&=&M \otimes_{A^{\rhd}} A^{\rhd}\left[\frac{1}{f_k}\right] 
\otimes_{A_0\left[\frac{f_0}{f_k},\cdots, \frac{f_n}{f_k}\right]}
\left(A_0\left[\frac{f_0}{f_k},\cdots, \frac{f_n}{f_k}\right]\right)^{Z}\\
&=&M \otimes_{A^{\rhd}} \left(1+I_0A_0\left[\frac{f_0}{f_k},\cdots, \frac{f_n}{f_k}\right]\right)^{-1} A^{\rhd}\left[\frac{1}{f_k}\right]\\
&=& (M \otimes \MO_A^{\Zar})(U_k).
\end{eqnarray*}
where we set 
{\small 
$$\left(A_0\left[\frac{f_0}{f_k},\cdots, \frac{f_n}{f_k}\right]\right)^{Z}
:=\left(1+I_0A_0\left[\frac{f_0}{f_k},\cdots, \frac{f_n}{f_k}\right]\right)^{-1}
\left(A_0\left[\frac{f_0}{f_k},\cdots, \frac{f_n}{f_k}\right]\right)$$
}
and the first equation follows from 
(\ref{e-sheafy2}), (\ref{e-sheafy3}) and 
the equation (\ref{e-FK-sheafy}) in Subsection~\ref{ss-FK-global}. 
By the same argument, we get 
$$\Gamma( \mathcal D_+(f_{k_1}f_{k_2}), \mathcal M
)= (M \otimes \MO_A^{\Zar})(U_{k_1k_2}).$$ 
Therefore, the sequence (\ref{e-sheafy}) coincides 
with the following sequence 
$$0 \to \Gamma(W, \mathcal M) 
\to \prod_{0 \leq k\leq n}  \Gamma(\mathcal D_+(f_k), \mathcal M) 
\to \prod_{0\leq k_1<k_2\leq n}  \Gamma(\mathcal D_+(f_{k_1}f_{k_2}), \mathcal M).$$
This is an exact suquence, since $\mathcal M=\iota^{-1}\sigma_*\widetilde M$ is a sheaf. 
Thus also (\ref{e-sheafy}) is exact. 
\end{proof}

\subsection{Zariskian adic spaces}\label{ss-zar-ad-sp}

Let $\mathcal V$ be the category of 
the triples $(X, \MO_X, (v_x)_{x \in X})$, 
where $X$ is a topological space, $\MO_X$ is a sheaf of topological rings, 
and each $v_x$ is a valuation of the stalk $\MO_{X, x}$. 
An arrow $f:X \to Y$ is a morphism in $\mathcal V$ 
if $f$ is a continuous map such that 
$f^{\sharp}:\MO_Y \to f_*\MO_X$ 
is a continuous ring homomorphism and that 
the valuation $v_{f(x)}$ is equivalent to the composition of 
$v_x$ and $f^{\sharp}_x:\MO_{Y, f(x)} \to \MO_{X, x}$.  

\begin{dfn}\label{d-zar-adic-sp}
For an affinoid ring $A$, 
the object 
$$(\Spa\,A, \MO_A^{\Zar}, (v_x)_{x \in \Spa\,A})$$ 
of $\mathcal V$ is called the {\em Zariskian adic space associated with} $A$. 
An {\em affinoid Zariskian adic space} is an object of $\mathcal V$ 
which is isomorphic to the Zariskian adic space associated with an affinoid ring. 

A {\em Zariskian adic space} is an object $(X, \MO_X, (v_x)_{x \in X})$ 
such that any point $x \in X$ has an open neighbourhood $U$ of $X$ 
to which the restriction $(U, \MO_X|_U, (v_x)_{x \in U})$ is 
an affinoid Zariskian adic space. 
A {\em morphism} of adic spaces is a morphism in $\mathcal V$. 
\end{dfn}

\subsection{Examples violating the Tate acylicity}\label{ss-cex-Tate}

\begin{nota}\label{n-cex}
Fix a prime number $p$ such that $p \neq 2$. 
In the following, we equip $\Q$ and $\Z_{(p)}$ 
with the $p$-adic topologies. 
We equip $\Q[t]$ the group topology induced by $\{p^n\Z_{(p)}[t]\}_{n \in \Z_{>0}}$. 
Then the pair
$$A:=(\Q[t], \Z_{(p)}[t])$$
is an affinoid ring.  
Let 
$$U_1:=R\left(\frac{1}{t+1}\right), \quad U_2:=R\left(\frac{1}{t-1}\right)$$
be rational subsets of $\Spa\,A$. 
We get 
$$U_1 \cap U_2=R\left(\frac{1}{t^2-1}\right)$$
We consider the following map 
\begin{eqnarray*}
\rho: \MO_A^{\Zar}(U_1) \times \MO_A^{\Zar}(U_2) &\to & 
\MO_A^{\Zar}(U_1 \cap U_2)\\
(f, g)&\mapsto &f|_{U_1 \cap U_2}-g|_{U_1 \cap U_2}.
\end{eqnarray*}
\end{nota}

\begin{lem}\label{l-cex-cover}
We use Notation~\ref{n-cex}. 
Then 
$$\Spa\,A=R\left(\frac{1}{t+1}\right) 
\cup R\left(\frac{1}{t-1}\right).$$ 
\end{lem}

\begin{proof}
Take $v \in \Spa(\Q[t], \Z_{(p)}[t]) \setminus R(\frac{1}{t+1})$. 
Then we have that $v(t+1) < 1=v(2)$, which implies 
$$v(t-1)=v(t+1-2)=1.$$
Thus we get $v \in R(\frac{1}{t-1})$, as desired. 
\end{proof}

\begin{prop}\label{p-non-surje}
We use Notation~\ref{n-cex}. 
Then $\rho$ is not surjective.
\end{prop}

\begin{proof}
The map $\rho$ can be written as 
$$\rho:\Q\left[t, \frac{1}{t+1}\right]^{\Zar} \times \Q\left[t, \frac{1}{t-1}\right]^{\Zar} 
\to \Q\left[t, \frac{1}{t^2-1}\right]^{\Zar}.$$
We embed all the rings appearing above into $\Q(t)$. 
We have that 
$$\frac{1}{1+\frac{p}{t^2-1}} \in 
\left(1+p\Z_{(p)}\left[t, \frac{1}{t^2-1}\right]\right)^{-1}
\Q\left[t, \frac{1}{t^2-1}\right]=\Q\left[t, \frac{1}{t^2-1}\right]^{\Zar}.$$
Assume that this element is in the image of $\rho$. 
Let us derive a contradiction. 
We can write 
$$ \frac{1}{1+\frac{p}{t^2-1}}=\frac{g_1(t, \frac{1}{t+1})}{1+pf_1(t, \frac{1}{t+1})}+
\frac{g_2(t, \frac{1}{t-1})}{1+pf_2(t, \frac{1}{t-1})}$$
for some $f_i(X, Y) \in \Z_{(p)}[X, Y]$ and $g_i(X, Y) \in \Q[X, Y]$. 
Taking the multiplication with the product of the denominators, 
we get 
{\small 
\begin{equation}\label{e-cex1}
\left(1+pf_1\middle(t, \frac{1}{t+1}\middle)\middle)\middle(1+pf_2\middle(t, \frac{1}{t-1}\middle)\right)=\left(1+\frac{p}{t^2-1}\right)g\left(t, \frac{1}{t^2-1}\right)
\end{equation}
}
for some $g(t, \frac{1}{t^2-1}) \in \Q[t, \frac{1}{t^2-1}]$. 
We have that $g(t, \frac{1}{t^2-1}) \in \Z_{(p)}[t, \frac{1}{t^2-1}]$. 
Indeed, otherwise we can find a positive integer $\nu$ such that 
$p^{\nu}g(t, \frac{1}{t^2-1}) \in \Z_{(p)}[t, \frac{1}{t^2-1}]$ and its modulo $p$ reduction 
is not zero, which contradicts the fact that $\F_p[t, \frac{1}{t^2-1}]$ 
is an integral domain, where $\F_p:=\Z/p\Z$.

\begin{claim}
$(1+\frac{p}{t^2-1})\Z_{(p)}[t, \frac{1}{t^2-1}]$ is a prime ideal of $\Z_{(p)}[t, \frac{1}{t^2-1}]$. 
\end{claim}

\begin{proof}(of Claim) 
We have that $t^2-1+p$ 
is an irreducible polynomial over $\Q$, which in turn implies that 
$$(t^2-1+p)\Z_{(p)}[t]$$
is a prime ideal of $\Z_{(p)}[t]$. 
In particular, we get $t^2-1 \not\in (t^2-1+p)\Z_{(p)}[t]$, 
since if $t^2-1 \in (t^2-1+p)\Z_{(p)}[t]$, then the residue ring 
$$\Z_{(p)}[t]/(t^2-1+p)\Z_{(p)}[t] \simeq \Z_{(p)}[t]/(t^2-1, p)\Z_{(p)}[t] \simeq \F_p[t]/(t^2-1)$$
is not an integral domain. 
For $S:=\{(t^2-1)^r\}_{r \geq 0}$, we have that 
{\small 
$$S^{-1}(\Z_{(p)}[t]/(t^2-1+p)\Z_{(p)}[t]) \simeq 
\Z_{(p)}\left[t, \frac{1}{t^2-1}\middle]\middle/(t^2-1+p)\Z_{(p)}\middle[t, \frac{1}{t^2-1}\right]$$
}
is an integral domain, as the image of $t^2-1$ in $\Z_{(p)}[t]/(t^2-1+p)\Z_{(p)}[t]$ is nonzero. 
Therefore, Claim holds. 
\end{proof}

Let us go back to the proof  of Proposition \ref{p-non-surje}. 
By Claim and (\ref{e-cex1}), one of 
$$1+pf_1\left(t, \frac{1}{t+1}\right)\quad \text{and}\quad 
1+pf_2\left(t, \frac{1}{t-1}\right)$$
is contained in $(1+\frac{p}{t^2-1})\Z_{(p)}[t, \frac{1}{t^2-1}]$. 
By symmetry, we may assume that the former case occurs. 
Then we can write 
\begin{equation}\label{e-cex2}
1+pf_1\left(t, \frac{1}{t+1}\right)=\left(1+\frac{p}{t^2-1}\right)
h\left(t, \frac{1}{t^2-1}\right)
\end{equation}
for some $h(X, Y) \in \Z_{(p)}[X, Y]$. 
For the additive $(t-1)$-adic valuation $v_{t-1}:\Q(t) \to \Z$ 
satisfying $v_{t-1}(t-1)=1$, 
it follows from (\ref{e-cex2}) that 
$v_{t-1}(h(t, \frac{1}{t^2-1}))\geq 1$, i.e. 
we can write 
\begin{equation}\label{e-cex3}
h\left(t, \frac{1}{t^2-1}\right)=(t-1)h_1\left(t, \frac{1}{t+1}\right)+...+(t-1)^kh_k\left(t, \frac{1}{t+1}\right)
\end{equation}
for some $h_1,\cdots, h_k \in \Z_{(p)}[t, \frac{1}{t+1}]$. 
Combining (\ref{e-cex2}) and (\ref{e-cex3}), we obtain 
\begin{equation}\label{e-cex4}
1+pf_1\left(t, \frac{1}{t+1}\right)=\left(t-1+\frac{p}{t+1}\right)
\widetilde{h}\left(t, \frac{1}{t+1}\right),
\end{equation}
for 
$$\widetilde{h}\left(t, \frac{1}{t+1}\right):=
h_1\left(t, \frac{1}{t+1}\right)+...+(t-1)^{k-1}h_k\left(t, \frac{1}{t+1}\right) 
\in \Z_{(p)}\left[t, \frac{1}{t+1}\right].$$
Substituting $t=1$ for (\ref{e-cex4}), 
we get an equation of rational numbers: 
\begin{equation}\label{e-cex5}
1+pf_1\left(1, \frac{1}{2}\right)=\frac{p}{2} \times 
\widetilde{h}\left(1, \frac{1}{2}\right).
\end{equation}
Since all of $\frac{1}{2}, f_1(1, \frac{1}{2})$ and 
$\widetilde{h}(1, \frac{1}{2})$ are contained in $\Z_{(p)}$, 
the $p$-adic valuations of the both hand sides of (\ref{e-cex5})
are different, which is absurd. 
This completes the proof of Proposition \ref{p-non-surje}.
\end{proof}
\begin{thm}\label{t-non-TA}
We use Notation~\ref{n-cex}. 
Then it holds that 
$$H^1(\Spa\,A, \MO_A^{\Zar}) \neq 0.$$ 
\end{thm}

\begin{proof}
Set $X:=\Spa\,A$. 
Let $j_1:U_1 \to X$, $j_2:U_2 \to X$ and $j_3:U_1 \cap U_2 \to X$ 
be the open immersions. 
Lemma~\ref{l-cex-cover} induces the following Mayer--Vietoris exact sequence: 
$$0 \to \MO^{\Zar}_A \to (j_1)_*(\MO^{\Zar}_A|_{U_1})
\times (j_2)_*(\MO^{\Zar}_A|_{U_2}) 
\to (j_3)_*(\MO^{\Zar}_A|_{U_1 \cap U_2}) \to 0.$$
It follows from Proposition~\ref{p-non-surje} that  
$H^1(X, \MO_A^{\Zar}) \neq 0$, as desired. 
\end{proof}

\section{Images to affine spectra}

In this section, we study the image of the natural map 
$$\theta:\Spa\,(A, A^+) \to \Spec\,A, \quad v \mapsto \Ker(v)$$
for an affinoid ring $(A, A^+)$. 
In Subsection~\ref{ss-characterise}, 
we prove that $A$ is Zariskian 
if and only if $\text{Im}(\theta)$ contains all the maximal ideals of $A$ (Theorem~\ref{t-characterise}). 
We conclude some criteria for faithful flatness 
(Corollary~\ref{c-ff-criterion}, Corollary~\ref{c-ff-cplt}). 
In Subsection~\ref{ss-non-surje}, 
we find a Zariskian affinoid ring $(A, A^+)$ such that 
$\theta$ is not surjective (Theorem \ref{t-non-surje}).

\subsection{Dominance onto maximal spectra}\label{ss-characterise}

The purpose of this subsection is to show Theorem~\ref{t-characterise}.

\begin{thm}\label{t-characterise}
Let $(A, A^+)$ be an affinoid ring. 
Then the following are equivalent. 
\begin{enumerate}
\item 
$A$ is Zariskian. 
\item 
An arbitrary maximal ideal of $A$ is contained in the image of 
the natural map 
$$\Spa\,(A, A^+) \to \Spec\,A, \quad v \mapsto \Ker(v).$$
\end{enumerate}
\end{thm}

\begin{proof}
Assume that (1) holds. 
Thanks to Lemma \ref{l-zar-nonempty}(2), 
we can find a point $v \in \Spa\,(A, A^+)$ such that $\m=\Ker(v)$. 
Hence, (2) holds.

Assume that (1) does not hold, i.e. $A$ is not Zariskian. 
We can find an element 
$x \in (1+A^{\circ\circ}) \setminus A^{\times}$. 
Then there exists a maximal ideal $\m$ containing $x$. 
Consider a commutative diagram of the induced maps: 
$$\begin{CD}
\Spa\,(A^{\Zar}, (A^+)^{\Zar}) @>g >> \Spec\,A^{\Zar}\\
@VV\beta V @VV \alpha V\\
\Spa\,(A, A^+) @>f >> \Spec\,A.
\end{CD}$$
Since the natural map $\beta$ is bijective by Theorem~\ref{t-comp-factor} and \cite[Proposition 3.9]{Hub93}, 
we have that 
$$\text{Im}(f)=\text{Im}(f \circ \beta)=\text{Im}(\alpha \circ g) 
\subset \text{Im}(\alpha).$$
By the choice of $\m$, we have that $\m \not\in \text{Im}(\alpha)$. 
Thus it holds that $\m \not\in \text{Im}(f)$, hence (2) does not hold. 
\end{proof}

\begin{cor}\label{c-ff-criterion}
Let $\varphi:(A, A^+) \to (B, B^+)$ be a continuous ring homomorphism 
of affinoid rings. 
Assume that $A$ is Zariskian and 
the induced map $\Spa\,(B, B^+) \to \Spa\,(A, A^+)$ is surjective. 
Then the following hold. 
\begin{enumerate}
\item 
Any maximal ideal of $A$ is contained in the image of 
the induced map $\Spec\,B \to \Spec\,A$. 
\item 
If $\varphi$ is flat, then $\varphi$ is faithfully flat. 
\end{enumerate}
\end{cor}

\begin{proof}
We have a natural commutative diagram: 
$$\begin{CD}
\Spa\,(B, B^+) @>\theta_B >> \Spec\,B\\
@VV\varphi^{\flat}V @VV\varphi^{\sharp}V\\
\Spa\,(A, A^+) @>\theta_A >> \Spec\,A.
\end{CD}$$
Therefore, (1) holds by Theorem~\ref{t-characterise}. 
We obtain the assertion (2) by (1) and \cite[Theorem 7.3(ii)]{Mat89}. 
\end{proof}

\begin{cor}\label{c-ff-cplt}
Let $\varphi:A \to B$ be a continuous ring homomorphism 
of f-adic rings. 
Assume that $A$ is Zariskian and 
the induced map $\widehat{\varphi}:\widehat{A} \to \widehat{B}$ 
is an isomorphism of topological rings. 
Then the following hold. 
\begin{enumerate}
\item 
Any maximal ideal of $A$ is contained in the image of 
the induced map $\Spec\,B \to \Spec\,A$. 
\item 
If $\varphi$ is flat, then $\varphi$ is faithfully flat. 
\end{enumerate}
\end{cor}

\begin{proof}
It follows from \cite[Proposition 3.9]{Hub93} that 
$\varphi^{\flat}:\Spa\,(B, B^{\circ}) \to \Spa\,(A, A^{\circ})$ is bijective. 
Therefore, the assertions follow from Corollary~\ref{c-ff-criterion}. 
\end{proof}

\begin{cor}
Let $(A, A^+)$ be an affinoid ring. 
Set $X:=\Spa\,(A, A^+)$. 
Let $X=\bigcup_{i \in I} U_i$ be 
a finite open cover where $U_i$ is a rational subset 
of $\Spa\,(A, A^+)$ for any $i \in I$. 
Then the induced ring homomorphism 
$$\rho:\MO_X^{\Zar}(X) \to \prod_{i \in I}\MO_X^{\Zar}(U_i)$$
is faithfully flat. 
\end{cor}

\begin{proof}
Since the induced map 
$$\rho^{\flat}:
\Spa\,\left(\prod_{i \in I}\MO_X^{\Zar}(U_i), \prod_{i \in I}\MO_X^{\Zar, +}(U_i)\right)
\to 
\Spa\,\left(\MO_X^{\Zar}(X), \MO_X^{\Zar, +}(X)\right)$$
is the same as the open cover $\coprod_{i \in I} U_i \to X$, 
we see that $\rho^{\flat}$ is surjective. 
Since $\rho$ is flat, it is faithfully flat by Corollary~\ref{c-ff-criterion}(ii). 
\end{proof}

\subsection{Non-surjective example}\label{ss-non-surje}

The purpose of this subsection is 
to show Theorem~\ref{t-non-surje}, 
which asserts that there exists a Zariskian f-adic ring $R$ such that 
the natural map 
$$\Spa\,(R, R^+) \to \Spec\,R,\quad v \mapsto \Ker(v)$$
is not surjective for any ring of integral elements $R^+$ of $R$. 
Our example is based on \cite[Remark (2) after Theorem 10.17]{AM69}, 
hence let us start by recalling its construction.

\begin{nota}\label{n-AM}
Let $A:=\{f:\R \to \R\,|\,f\text{ is of class }C^{\infty}\}.$ 
For 
$$\varphi:A \to \R,\quad f \mapsto f(0),$$
we set $\m:=\Ker(\varphi)$. 
Since $\varphi$ is a surjective ring homomorphism to a field $\R$, 
it holds that $\m$ is a maximal ideal of $A$. 
By Taylor's theorem, we get $\m=xA$. 
It follows again from Taylor's theorem that 
$$\bigcap_{n=1}^{\infty} \m^n=\{f \in A\,|\,f(0)=f'(0)=f''(0)=\cdots=0\}.$$
We equip $A$ with the $\m$-adic topology. 
Since $\m$ is a finitely generated ideal, we have that $A$ is an f-adic ring. 
Recall that $A^{\Zar}=(1+\m)^{-1}A$ and 
let $\alpha:A \to A^{\Zar}$ be the natural continuous ring homomorphism, 
where the topology on $A^{\Zar}$ coincides with the $\m A^{\Zar}$-adic topology. 
We can directly check that 
$$\Ker(\alpha)=\{f \in A\,|\, f|_U=0\text{ for some open neighbourhood } U\text{ of } 0 \in \R\}.$$
In particular, it follows that 
$e^{-1/x^2} \in \left(\bigcap_{n=1}^{\infty} \m^n\right) \setminus \Ker(\alpha)$ 
and 
$$\alpha(e^{-1/x^2}) \in \left(\bigcap_{n=1}^{\infty} \m^n A^{\Zar}\right) \setminus \{0\}.$$
\end{nota}

\begin{lem}\label{l-non-surje}
We use Notation~\ref{n-AM}. 
Then the following hold. 
\begin{enumerate}
\item $A$ is reduced. 
\item $A^{\Zar}$ is reduced. 
\item There exists a prime ideal $\p$ of $A^{\Zar}$ 
such that 
$$\bigcap_{n=1}^{\infty} \m^n A^{\Zar} \not\subset \p.$$ 
\item The closed immersion 
$$\pi^{\sharp}:\Spec\,((A^{\Zar})^{\hd}) \to \Spec\,A^{\Zar}$$
is not surjective, where $\pi^{\sharp}$ is the morphism 
induced by the natural surjective ring homomorphism 
$\pi:A^{\Zar} \to (A^{\Zar})^{\hd}$.
\end{enumerate}
\end{lem}

\begin{proof}
We first show (1). 
Take $f \in A$ and assume that $f^n=0$ for some positive integer $n$. 
Since $\R$ is reduced, we have that 
the equation $f(a)^n=0$ implies that $f(a)=0$ for any $a \in \R$. 
Therefore, we get $f=0$. 
Thus (1) holds. 
The assertion (2) follows from (1) (cf. \cite[Corollary 3.12]{AM69}).

We now show (3). 
It follows from (2) and $\bigcap_{n=1}^{\infty} \m^n A^{\Zar} \neq 0$ that 
$$\bigcap_{n=1}^{\infty} \m^n A^{\Zar} \not\subset \{0\}=\sqrt{0}=\bigcap_{\p \in \Spec\,A^{\Zar}} \p.$$
In particular, there exists a prime ideal $\p$ of $A^{\Zar}$ 
such that $\bigcap_{n=1}^{\infty} \m^n A^{\Zar} \not\subset \p$. 
Thus (3) holds.

Let us show (4). 
Since $(A^{\Zar})^{\hd}=A^{\Zar}/\left( \bigcap_{n=1}^{\infty} \m^n A^{\Zar}\right)$, 
the image of $\pi^{\sharp}$ consists of the prime ideals of $A^{\Zar}$ 
containing $\bigcap_{n=1}^{\infty} \m^n A^{\Zar}$. 
Thus (4) follows from (3). 
\end{proof}

\begin{thm}\label{t-non-surje}
There exists a Zariskian f-adic ring $R$ such that the natural map 
$$\theta:\Spa\,(R, R^+) \to \Spec\,R, \quad v \mapsto \Ker(v)$$
is not surjective for any ring of integral elements $R^+$ of $R$. 
\end{thm}

\begin{proof}
We use Notation~\ref{n-AM}. 
Set $R:=A^{\Zar}$. 
Fix an arbitrary ring of integral elements $R^+$ of $R$. 
We obtain a commutative diagram of natural maps: 
$$\begin{CD}
\Spa\,(R^{\hd}, (R^+)^{\hd}) @>\theta^{\hd}>> \Spec\,R^{\hd}\\
@VV\pi^{\flat}V @VV\pi^{\sharp} V\\
\Spa\,(R, R^+) @>\theta >> \Spec\,R.
\end{CD}$$
It follows from \cite[Proposition 3.9]{Hub93} that $\pi^{\flat}$ is bijective. 
Since $\pi^{\sharp}$ is not surjective by Lemma~\ref{l-non-surje}, 
neither is $\theta$. 
\end{proof}

\end{document}